\documentclass[graybox]{svmult}

\smartqed
\usepackage{mathptmx}       % selects Times Roman as basic font
\usepackage{helvet}         % selects Helvetica as sans-serif font
\usepackage{courier}        % selects Courier as typewriter font
\usepackage{type1cm}        % activate if the above 3 fonts are
\usepackage{makeidx}         % allows index generation
\usepackage{graphicx}        % standard LaTeX graphics tool
\usepackage{multicol}        % used for the two-column index
\usepackage[bottom]{footmisc}% places footnotes at page bottom

\usepackage{amsmath,graphicx}
\usepackage{graphicx}
\usepackage{epsfig}
\usepackage{sidecap}
\usepackage{amsmath}
\usepackage{amssymb}
\usepackage{subfigure}
\usepackage{wrapfig}
\usepackage{times,color}
\usepackage{cite,footnote,xspace,syntonly,theorem,algorithm,algorithmic,bm}

\input{mysymbol.sty}
\spnewtheorem{assumption}{Assumption}{\bf}{\it}

\makeindex             % used for the subject index
\newcommand{\calV}{{\cal V}}

\newcommand{\calN}{{\cal N}}

\newcommand{\calZ}{{\cal Z}}

\newcommand{\calS}{{\cal S}}

\def \cP {\mathcal{P}}

\def \bY {{\bf Y}}
\def \bR {{\bf R}}
\def \bS {{\bf S}}
\def \bX {{\bf X}}
\def \bA {{\bf A}}

\def \bI {{\bf I}}

\def \bE {{\bf E}}

\def \bZ {{\bf Z}}

\def \bH {{\bf H}}

\def \bQ {{\bf Q}}

\def \bP {{\bf P}}

\def \bv {{\bf v}}
\def \bx {{\bf x}}
\def \by {{\bf y}}

\def \br {{\bf r}}

\def \bv {{\bf v}}
\def \bs {{\bf s}}

\def \rank {\textsf{rank}}
\def \Tr {\textsf{Tr}}

\begin{document}

\title*{Decentralized Learning for Wireless Communications and Networking}
\author{Georgios B. Giannakis, Qing Ling,  Gonzalo Mateos, Ioannis D. Schizas\\
and Hao Zhu}
\institute{Georgios B. Giannakis \at University of Minnesota, 117 Pleasant Str.,
Minneapolis, MN 55455, \email{georgios@umn.edu}
    \and Qing Ling \at University of Science and Technology of China, 443 Huangshan Road, Hefei, Anhui, China 230027, \email{qingling@mail.ustc.edu.cn}
    \and Gonzalo Mateos \at University of Rochester, 413 Hopeman Engineering
Building, Rochester, NY 14627, \email{gmateosb@ece.rochester.edu}
    \and Ioannis D. Schizas \at University of Texas at Arlington, 416 Yates Street,
Arlington, TX 76011, \email{schizas@uta.edu}
    \and Hao Zhu \at University of Illinois at Urbana-Champaign, 4058 ECE Building, 306 N. Wright
Street, Urbana, IL 61801, \email{haozhu@illinois.edu}}
\maketitle

\abstract{This chapter deals with decentralized learning algorithms for
in-network
processing of graph-valued data. A generic learning problem is formulated and
recast
into a separable form, which is iteratively minimized using the
alternating-direction method of multipliers (ADMM) so as to gain the
desired degree of parallelization. Without exchanging elements from the
distributed training sets and keeping inter-node communications at affordable
levels, the local (per-node) learners consent to the desired quantity inferred
globally, meaning the one obtained if the entire training data set were centrally
available. Impact of the decentralized learning framework to contemporary wireless
communications and networking tasks is illustrated through case studies including
target tracking using wireless sensor networks, unveiling Internet traffic anomalies,
power system state estimation, as well as spectrum cartography for wireless
cognitive radio networks.}

\section{Introduction}\label{sec:Introduction}

This chapter puts forth an optimization framework for learning
over networks, that entails decentralized processing of training data
acquired by interconnected nodes.
Such an approach is of paramount importance when communication of
training data to a central processing unit
is prohibited due to e.g., communication cost or privacy reasons. The so-termed
in-network
processing paradigm for decentralized
learning is based on successive refinements of local model parameter estimates
maintained at individual network nodes. In a nutshell, each iteration of this
broad class of fully
decentralized algorithms comprises: (i) a communication step where nodes
exchange information
with their neighbors through e.g., the shared wireless medium or Internet
backbone; and (ii) an
update step where each node uses this information to refine its local estimate.
Devoid of hierarchy and with their decentralized in-network processing,
local e.g., estimators should eventually consent to the global estimator
sought, while fully exploiting existing spatiotemporal correlations to
maximize estimation performance. In most cases, consensus can formally
be attained asymptotically in time. However, a finite number of iterations will
suffice to
obtain results that are sufficiently accurate for all practical purposes.

In this context, the approach followed here entails reformulating a generic
learning task as a convex
constrained optimization problem, whose structure lends itself
naturally to decentralized implementation over a network graph. It is then
possible to
capitalize upon this favorable structure by resorting to the
alternating-direction method of multipliers (ADMM), an
iterative optimization method that can be traced back
to~\cite{Glowinski_Marrocco_ADMM_1975}
(see also~\cite{Gabay_Mercier_ADMM_1976}),
and which is specially well-suited for parallel
processing~\cite{bertsi97book,Boyd_ADMM}. This way simple decentralized recursions
become available to update each node's local estimate, as well as
a vector of dual prices through which network-wide agreement is effected.

\runinhead{Problem statement.} Consider a network of $n$ nodes in which scarcity
of power and bandwidth resources encourages only single-hop inter-node
communications, such that the $i$-th node communicates solely
with nodes $j$ in its single-hop neighborhood ${\cal N}_i$. Inter-node links
are assumed symmetric, and the network is modeled
as an undirected graph whose vertices are the nodes and its
edges represent the available communication links. As it will become clear
through
the different application domains studied here, nodes could be wireless
sensors, wireless access points (APs), electrical buses, sensing cognitive radios,
or routers, to name a few examples.
Node $i$ acquires $m_i$ measurements stacked
in the vector $\mathbf{y}_i\in\mathbb{R}^{m_i}$ containing information about
the unknown model parameters in $\mathbf{s}\in\mathbb{R}^{p}$,
which the nodes need to estimate. Let
$\mathbf{y}:=[\mathbf{y}_1^\top,\ldots,\mathbf{y}_n^\top]^\top\in\mathbb{R}^{\sum_{i}m_i}$
collect measurements acquired across the entire network. Many popular centralized
schemes obtain an estimate $\hat{\mathbf{s}}$ as follows
\begin{equation}\label{Eq:Centr_Est}
\hat{\mathbf{s}}\in\arg\min_{\mathbf{s}}\textstyle\sum_{i=1}^{n}f_{i}(\mathbf{s};\mathbf{y}_i).
\end{equation}
In the decentralized learning problem studied here though, the summands $f_{i}$
are assumed to be local cost functions only known to each node $i$. Otherwise
sharing this information with a centralized processor, also referred to as fusion center
(FC), can be challenging in various applications of interest, or, it may be
even impossible in e.g., wireless sensor networks (WSNs) operating under
stringent power budget constraints. In other cases such as the Internet or
collaborative healthcare studies, agents may not be willing to
share their private training data $\by_i$ but only the learning results.
Performing the optimization \eqref{Eq:Centr_Est} in a centralized fashion raises robustness
concerns as well, since the central processor represents an isolated point of
failure.

In this context, the objective of this chapter is to develop a decentralized
algorithmic framework for learning tasks, based on in-network processing of
the locally available data.
The described setup naturally suggests three characteristics that the algorithms
should exhibit: c1) each node $i=1,\ldots,n$ should obtain an estimate of
$\bs$, which coincides with the corresponding solution $\hat{\bs}$
of the centralized estimator \eqref{Eq:Centr_Est} that uses the entire data
$\{\by_i\}_{i=1}^n$; c2) processing per node should be kept as
simple as possible; and c3) the overhead for inter-node
communications should be affordable and confined to single-hop
neighborhoods. It will be argued that such an ADMM-based algorithmic framework
can
be useful for contemporary applications in the domain of wireless communications
and networking.

{\runinhead{Prior art.} Existing decentralized solvers of
\eqref{Eq:Centr_Est} can be classified in two categories:
C1) those obtained by modifying centralized algorithms and
operating in the primal domain; and C2) those handling an
equivalent constrained form of \eqref{Eq:Centr_Est} (see
\eqref{Eq:Separanle_Estimation} in Section \ref{sec:ADMM}),
and operating in the primal-dual domain.

Primal-domain algorithms under C1 include the (sub)gradient method
and its variants \cite{Nedic2009, Ram2010, Yuan2013, Jakovetic2013},
the incremental gradient method \cite{Rabbat2005-inc}, the proximal
gradient method \cite{Chen2012}, and the dual averaging method
\cite{Duchi2012, Tsianos2012-acc}. Each node in these methods,
averages its local iterate with those of neighbors and
descends along its local negative (sub)gradient direction.
However, the resultant algorithms are limited to inexact convergence when
using constant stepsizes \cite{Nedic2009,Yuan2013}. If diminishing
stepsizes are employed instead, the algorithms can achieve exact
convergence at the price of slowing down speed
\cite{Jakovetic2013,Rabbat2005-inc,Duchi2012}. A constant-stepsize
exact first-order algorithm is also available to achieve fast and
exact convergence, by correcting error terms in the distributed
gradient iteration with two-step historic information
\cite{Shi2014-extra}.

Primal-dual domain algorithms under C2 solve an equivalent
constrained form of \eqref{Eq:Centr_Est}, and thus drive local
solutions to reach global optimality. The dual decomposition method
is hence applicable because (sub)gradients of the dual function depend
on local and neighboring iterates only, and can thus be computed
without global cooperation \cite{Rabbat2005}. ADMM modifies the dual
decomposition by regularizing the constraints with a quadratic term,
which improves numerical stability as well as rate of convergence,
as will be demonstrated later in this chapter. Per ADMM iteration,
each node solves a subproblem that can be demanding. Fortunately,
these subproblems can be solved inexactly by running one-step
gradient or proximal gradient descent iterations, which
markedly mitigate the computation burden
\cite{Ling2014-icassp,Chang2014-icassp}. A sequential distributed ADMM algorithm
can be found in~\cite{Wei2012}.

\runinhead{Chapter outline.} The remainder of this chapter is organized as
follows.
Section \ref{sec:ADMM} describes a generic ADMM framework for decentralized
learning
over networks, which is at the heart of all algorithms described in the chapter
and was
pioneered in~\cite{sg06asilomar,srg08tsp} for in-network estimation using WSNs. Section
\ref{sec:Batch}  focuses on batch estimation as well as (un)supervised
inference,
while Section \ref{sec:Adaptive} deals with decentralized adaptive estimation
and tracking schemes where network nodes collect data sequentially in time.
Internet traffic anomaly detection and spectrum cartography for
wireless CR networks serve as motivating applications for the
sparsity-regularized
rank minimization algorithms developed in Section \ref{sec:Sparse}. Fundamental
results
on the convergence and convergence rate of decentralized ADMM are stated in
Section \ref{sec:Convergence}.

\section{In-Network Learning with ADMM in a Nutshell}\label{sec:ADMM}

Since local summands in \eqref{Eq:Centr_Est} are coupled through a \emph{global}
variable
$\mathbf{s}$, it is not straightforward to decompose the
unconstrained optimization problem in \eqref{Eq:Centr_Est}. To overcome this
hurdle,
the key idea is to introduce local variables $\calS:=\{\mathbf{s}_i\}_{i=1}^n$
which
represent local estimates of $\mathbf{s}$ per network node $i$~\cite{sg06asilomar,srg08tsp}.
Accordingly, one can
formulate the \emph{constrained} minimization problem
\begin{equation}\label{Eq:Separanle_Estimation}
     \{\hat{\mathbf{s}}_i\}_{i=1}^{n}\in\arg\min_{\calS}\textstyle\sum_{i=1}^{n}
    f_{i}(\mathbf{s}_i;\mathbf{y}_i), \;\;\;\; \textrm{s. to }~~~
    \mathbf{s}_i=\mathbf{s}_j,\;j\in{\cal N}_i.
\end{equation}
The ``consensus'' equality constraints in \eqref{Eq:Separanle_Estimation} ensure
that
local estimates coincide within neighborhoods. Further, if the graph is
connected then
consensus naturally extends to the whole network,
and  it turns out that  problems \eqref{Eq:Centr_Est} and
\eqref{Eq:Separanle_Estimation}
are equivalent in the sense that
$\hat{\mathbf{s}}=\hat{\mathbf{s}}_1=\ldots=\hat{\mathbf{s}}_{n}$~\cite{srg08tsp}.
Interestingly, the formulation in \eqref{Eq:Separanle_Estimation} exhibits a
separable structure that is
amenable to decentralized minimization. To leverage this favorable structure,
the alternating direction method of multipliers (ADMM), see e.g., \cite[pg.
253-261]{bertsi97book}, can be employed here to minimize
\eqref{Eq:Separanle_Estimation}
in a decentralized fashion. This procedure will yield a
distributed estimation algorithm whereby local iterates
$\mathbf{s}_i(k)$, with $k$ denoting iterations, provably converge
to the centralized estimate $\hat{\mathbf{s}}$ in \eqref{Eq:Centr_Est};
see also Section \ref{sec:Convergence}.

To facilitate application of ADMM, consider the auxiliary variables
$\calZ:=\{\mathbf{z}_i^{j}\}_{j\in{\cal N}_i}$, and reparameterize the
constraints in \eqref{Eq:Separanle_Estimation} with the equivalent ones
\begin{align}\label{Eq:Constr_Equi}
     \{\hat{\mathbf{s}}_i\}_{i=1}^{n}&{}\in\arg\min_{\calS}\textstyle\sum_{i=1}^{n}
    f_{i}(\mathbf{s}_i;\mathbf{y}_i),\nonumber\\
    \textrm{s. to }&~~~\mathbf{s}_i=\mathbf{z}_i^{j}\textrm{
        and   }\mathbf{s}_{j}=\mathbf{z}_i^{j},\;\;
    i=1,\ldots,n,\;\;j\in{\cal  N}_i,\;\;i\neq j.
\end{align}
Variables $\mathbf{z}_i^{j}$ are only used to derive the local
recursions but will be eventually eliminated.  Attaching Lagrange multipliers
$\calV:=\{\{\bar{\mathbf{v}}_{i}^j\}_{j\in\calN_i},\{\tilde{\mathbf{v}}_i^j\}_{j\in\calN_i}\}_{i=1}^n$
to the constraints \eqref{Eq:Constr_Equi},
consider the augmented Lagrangian function
\begin{align}\label{Eq:Lagrangian_Function}
    L_c[\calS,\calZ,\calV] = \sum_{i=1}^{n}f_{i}(\mathbf{s}_i;\mathbf{y}_i)
    &+\sum_{i=1}^{n}\sum_{j\in{\cal N}_{i}}\left[(\bar{\mathbf{v}}_{i}^{j})^\top
    (\mathbf{s}_{i}-\mathbf{z}_{i}^{j})+(\tilde{\bv}_i^{j})^\top
    (\mathbf{s}_{j}-\mathbf{z}_{i}^{j})\right]\nonumber\\
    %Qudratic penalty coefficients associated  to the equality constraints
    &+\frac{c}{2}\sum_{i=1}^{n}\sum_{j\in{\cal N}_i}
\left[\|\mathbf{s}_{i}-\mathbf{z}_{i}^{j}\|^{2}+\|\mathbf{s}_j-\mathbf{z}_{i}^{j}\|^2\right]
\end{align}
where the constant $c>0$ is a penalty coefficient.
To minimize \eqref{Eq:Separanle_Estimation}, ADMM entails an iterative procedure
comprising three steps per iteration
$k=1,2,\ldots$
\begin{description}
    \item [{\bf [S1]}] \textbf{Multiplier updates:}\begin{align}
         \hspace{-1cm}\bar{\bbv}_i^{j}(k)&=\bar{\bbv}_{i}^{j}(k-1)+
        c[\bbs_i(k)-\bbz_i^{j}(k)] \nonumber\\
         \hspace{-1cm}\tilde{\bbv}_i^{j}(k)&=\tilde{\bbv}_{i}^{j}(k-1)+
        c[\bbs_j(k)-\bbz_i^{j}(k)].\nonumber
    \end{align}

    \item [{\bf [S2]}]  \textbf{Local estimate updates:}
    \begin{equation}\hspace{-1cm}\calS(k+1)=
         \mbox{arg}\:\min_{\calS}L_c\left[\calS,\calZ(k),\calV(k)\right].
        \nonumber\end{equation}

    \item [{\bf [S3]}] \textbf{Auxiliary variable updates:}
    \begin{equation}\hspace{-1cm}\calZ(k+1)=
         \mbox{arg}\:\min_{\calZ}L_{c}\left[\calS(k+1),\calZ,\calV(k)\right]
        \nonumber\end{equation}
\end{description}
where $i=1,\ldots,n$ and $j\in\calN_{i}$ in [S1]. Reformulating the generic
learning problem \eqref{Eq:Centr_Est} as \eqref{Eq:Constr_Equi} renders
the augmented Lagrangian in \eqref{Eq:Lagrangian_Function} highly
decomposable. The separability comes in two flavors, both
with respect to the sets $\calS$ and $\calZ$ of primal variables, as
well as across nodes $i=1,\ldots,n$. This in turn leads to highly
parallelized, simplified recursions corresponding to the
aforementioned steps [S1]-[S3]. Specifically, as detailed in
e.g.,~\cite{srg08tsp,sgrr08tsp,SMG_D_LMS,pfacgg10jmlr,mateos_dlasso,mmg13tsp},
it follows that if the multipliers are initialized to zero, the
ADMM-based decentralized algorithm reduces to the following updates carried
out locally at every node
\begin{svgraybox}
    \noindent \textbf{In-network learning algorithm at node $i$, for
$k=1,2,\ldots$:}
    \begin{align}
\label{Eq:vupdate}\bv_{i}(k)&=\bv_{i}(k-1)+c\sum_{j\in\calN_i}[\bs_i(k)-\bs_{j}(k)]\\
        \label{Eq:supdate}\bs_i(k+1)&=\arg\min_{\bs_i}\left\{
        f_{i}(\mathbf{s}_i;\mathbf{y}_i)+\bv_i^{\top}(k)\bs_i
        +c\sum_{j\in\calN_i}
         \left\|\bs_i-\frac{\bs_i(k)+\bs_{j}(k)}{2}\right\|^{2}\right\}
    \end{align}
\end{svgraybox}
\noindent where
$\bv_i(k):=2\sum_{j\in\calN_i}\bar{\bv}_{i}^{j}(k)$,
and all initial values are set to zero.

Recursions \eqref{Eq:vupdate} and \eqref{Eq:supdate} entail local
updates, which comprise the general purpose ADMM-based decentralized
learning algorithm. The inherently redundant set
of auxiliary variables in $\calZ$ and corresponding multipliers
have been eliminated. Each node, say the $i$-th one, does not need
to \textit{separately} keep track of all its non-redundant multipliers
$\{\bar{\bbv}_i^{j}(k)\}_{j\in\calN_i}$, but only to update
the (scaled) sum $\bbv_i(k)$. In the end, node $i$ has to
store and update only two $p$-dimensional vectors, namely
$\{\bs_i(k)\}$ and $\{\bv_i(k)\}$. A unique feature of in-network processing is
that nodes communicate their updated
local estimates $\{\bs_i\}$ (and not their
raw data $\by_i$) with their neighbors, in order
to carry out the tasks \eqref{Eq:vupdate}-\eqref{Eq:supdate} for the next
iteration.

As elaborated in Section \ref{sec:Convergence}, under mild assumptions
on the local costs one can establish that
$\lim_{k\to\infty}\bs_i(k)=\hat{\bs}$, for
$i=1,\ldots,n$. As a result, the algorithm asymptotically attains consensus and
the performance of the centralized
estimator [cf. \eqref{Eq:Centr_Est}].

%%%%%%%%%%%%%%%%%%%%%%%%%%%%%%%%%%%%%%%%%%%%%%%%%%%%%%%%%%%%%%%%%%%%%%
%                                                                    %
%               Subsection: Hao's part                               %
%  - Classification (Linear-Bayes, SVM, DiVA, Clustering)            %
%                 4 pages                                            %
%%%%%%%%%%%%%%%%%%%%%%%%%%%%%%%%%%%%%%%%%%%%%%%%%%%%%%%%%%%%%%%%%%%%%%

\section{Batch In-Network Estimation and Inference}\label{sec:Batch}

\subsection{Decentralized Signal Parameter Estimation}\label{subsec:Estimation}

%Traditional estimation techniques are developed under the critical assumption
%that the measurements $\mathbf{x}$
%are available at a centralized location. This fusion center is also a central
%processor which runs
%specialized algorithms to estimate unknown parameters of interest, namely
%$\mathbf{s}\in\mathbb{R}^{r\times 1}$;
%the latter vector could correspond to e.g., a thermal source, or a target
%position.
Many workhorse estimation schemes such as maximum likelihood estimation (MLE),
least-squares estimation (LSE), best linear unbiased estimation (BLUE),
as well as linear minimum mean-square error estimation (LMMSE) and
the maximum a posteriori (MAP) estimation, all can be formulated as a minimization
task similar
to \eqref{Eq:Centr_Est}; see e.g.~\cite{Estimation_Theory}.
However, the corresponding centralized estimation algorithms fall short in
settings where
both the acquired measurements and computational capabilities are distributed
among multiple
spatially scattered sensing nodes, which is the case with WSNs.
Here we outline a novel batch decentralized optimization framework building on
the ideas in Section~\ref{sec:ADMM},
that formulates the desired estimator as the solution of a separable constrained
convex minimization problem
tackled via ADMM; see e.g., \cite{bertsi97book,Boyd_ADMM,srg08tsp,sgrr08tsp} for
further details on the algorithms outlined
here.

Depending on the estimation technique utilized, the local cost functions
$f_i(\cdot)$ in \eqref{Eq:Centr_Est} should be chosen
accordingly, see e.g., \cite{Estimation_Theory,srg08tsp,sgrr08tsp}. For
instance, when $\mathbf{s}$ is assumed to be an unknown
deterministic vector, then:
\begin{itemize}
    \item If $\hat{\mathbf{s}}$ corresponds to the centralized MLE then
     $f_i(\mathbf{s};\mathbf{y}_i)=-\ln[p_i(\mathbf{y}_i;\mathbf{s})]$ is the
negative log-likelihood capturing the data probability
    density function (pdf), while the network-wide data
$\{\mathbf{y}_i\}_{i=1}^n$ are assumed statistically independent.

    \item If $\hat{\mathbf{s}}$ corresponds to the BLUE (or weighted least-squares
estimator) then
$f_i(\mathbf{s};\mathbf{y}_i)=(1/2)\|\bm{\Sigma}_{y_i}^{-1/2}(\mathbf{y}_i-\mathbf{H}_{i}\mathbf{s})\|^2$,  where $\bm{\Sigma}_{y_i}$ denotes the covariance of the data
$\mathbf{y}_i$, and $\mathbf{H}_i$ is a known fitting matrix.
\end{itemize}
When $\mathbf{s}$ is treated as a random vector, then:
\begin{itemize}
    \item If $\hat{\mathbf{s}}$ corresponds to the centralized MAP estimator then
$f_i(\mathbf{s};\mathbf{y}_i)=-(\ln[p_i(\mathbf{y}_i|\mathbf{s})]+n^{-1}\ln[p(\mathbf{s})])$
accounts for the data pdf, and $p(\mathbf{s})$ for the prior pdf of $\mathbf{s}$, while
data $\{\mathbf{y}_i\}_{i=1}^{n}$ are assumed conditionally independent given $\mathbf{s}$.

    \item If $\hat{\mathbf{s}}$ corresponds to the centralized LMMSE then
$f_i(\mathbf{s};\mathbf{y}_i)=(1/2)\|\mathbf{s}-n\boldsymbol{\Sigma}_{sy_i}\mathbf{u}^i\|_2^2$,
where $\bm{\Sigma}_{sy_i}$ denotes the cross-covariance of $\mathbf{s}$ with
$\mathbf{y}_i$, while $\mathbf{u}^i$ stands for the $i$-th $m_i\times 1$ block
subvector of $\mathbf{u}=\bm{\Sigma}_{y}^{-1}\mathbf{y}$.
\end{itemize}
Substituting in \eqref{Eq:supdate} the specific $f_i(\mathbf{s};\mathbf{y}_i)$
for each of the
aforementioned estimation tasks, yields a family of batch ADMM-based
decentralized estimation algorithms. The decentralized BLUE algorithm will be
described in this section as an example of decentralized linear estimation.

Recent advances in cyber-physical systems have also stressed the need for
decentralized nonlinear least-squares (LS) estimation. Monitoring the power grid for
instance, is challenged by the nonconvexity arising from the nonlinear AC
power flow model; see e.g., \cite[Ch. 4]{Wollenberg-book}, while the
interconnection across local transmission systems motivates their
operators to collaboratively monitor the global system state.
Interestingly, this nonlinear (specifically quadratic) estimation task
can be convexified to a semidefinite program (SDP) \cite[pg. 168]{Boyd_Convex},
for which a decentralized semidefinite programming (SDP) algorithm can
be developed by leveraging the batch ADMM; see also~\cite{Wen2010} for an ADMM-based
	centralized SDP precursor.

\subsubsection{Decentralized BLUE}
\label{subsubsec:DBLUE}

The minimization involved in \eqref{Eq:supdate} can be performed locally at
sensor $i$ by employing numerical
optimization techniques~\cite{Boyd_Convex}.
There are cases where the minimization in \eqref{Eq:supdate} yields a closed-form
and easy to implement
updating formula for $\mathbf{s}_i(k+1)$. If for example network nodes wish
to find the BLUE estimator in a distributed fashion, the local cost is
$f_i(\mathbf{s};\mathbf{y}_i)=(1/2)\|\bm{\Sigma}_{y_i}^{-1/2}
(\mathbf{y}_i-\mathbf{H}_{i}\mathbf{s})\|^2$, and \eqref{Eq:supdate}
becomes a strictly convex unconstrained quadratic program which admits
the following closed-form solution (see details in \cite{srg08tsp,MSG_D_RLS})
\begin{align}\label{Eq:D_BLUE}
\mathbf{s}_{i}(k+1)&=\left(\mathbf{H}_i^\top\bm{\Sigma}_{y_i}^{-1}\mathbf{H}_i+2c|{\cal
N}_i|\mathbf{I}_{p}\right)^{-1}\left[
\mathbf{H}_i^\top\bm{\Sigma}_{y_i}^{-1}\mathbf{y}_i-\mathbf{v}_i(k)+c\sum_{j\in{\cal
N}_i}   \left(\mathbf{s}_i(k)+\mathbf{s}_j(k)\right)\right].
\end{align}
The pair \eqref{Eq:vupdate} and  \eqref{Eq:D_BLUE} comprise the decentralized
(D-) BLUE algorithm~\cite{sg06asilomar,srg08tsp}.
For the special case where each node acquires unit-variance scalar observations
$y_i$, there is no fitting matrix
and $s$ is scalar (i.e., $p=1$); D-BLUE offers a decentralized algorithm to
obtain the network-wide sample average $\hat{s}=(1/n)\sum_{i=1}^n y_i$. The update
rule for the local estimate is obtained by suitably specializing \eqref{Eq:D_BLUE} to
\begin{align}\label{Eq:D_BLUE_averaging}
    s_{i}(k+1)&=\left(1+2c|{\cal N}_i|\right)^{-1}\left[
y_i-v_i(k)+c\sum_{j\in{\cal N}_i}
    \left(s_i(k)+s_j(k)\right)\right].
\end{align}
Different from existing distributed averaging approaches~\cite{Barbarossa_Scu_Coupled_Osci,dimakis10,Consensus_Averaging,xbk07jpdc},
the ADMM-based one originally proposed in~\cite{sg06asilomar,srg08tsp} allows the
decentralized computation of general nonlinear estimators that may be not
available in closed form and cannot be expressed as ``averages.'' Further,
the obtained recursions exhibit robustness in the presence of additive noise in
the inter-node communication links.

\subsubsection{Decentralized SDP}
\label{subsubsec:DSDP}

Consider now that each scalar $y_i^\ell$ in $\by_i$ adheres to a quadratic measurement model in $\bbs$ plus additive Gaussian noise, where the centralized MLE requires solving a nonlinear least-squares problem. To tackle the nonconvexity due to the quadratic dependence, the task of estimating the state  $\bs$ can be reformulated as that of estimating the outer-product matrix $\bbS :=\bs \bs^\top$. In this reformulation $y_i^\ell$ is a linear function of $\bS$, given by $\Tr(\bH_i^\ell\bS)$ with a known matrix $\bH_i^\ell$~\cite{hzgg_jstsp14}. Motivated by the separable structure in \eqref{Eq:Constr_Equi}, the nonlinear estimation problem can be similarly formulated as
\begin{align}\label{Eq:cen_SDP}
     \{\hat{\mathbf{S}}_i\}_{i=1}^{n}&{}\in\arg\min \sum_{i=1}^{n}
     \sum_\ell \left[y_i^\ell - \Tr(\bH_i^\ell\bS) \right]^2,\nonumber\\
    \textrm{s. to }&~~~\mathbf{S}_i=\mathbf{Z}_i^{j}\textrm{
        and   }\mathbf{S}_{j}=\mathbf{Z}_i^{j},\;\;
    i=1,\ldots,n,\;\;j\in{\cal  N}_i,\;\;i\neq j \nonumber\\
    &~~~ \bS_i \succeq \mathbf{0}\textrm{
            and   }\rank(\bS_i)=1,\;\;
                i=1,\ldots,n~
\end{align}
where the positive-semidefiniteness and rank constraints ensure that each matrix $\bS_i$ is an outer-product matrix. By dropping the non-convex rank constraints, the problem \eqref{Eq:cen_SDP} becomes a convex semidefinite program (SDP), which can be solved in a decentralized fashion by adopting the batch ADMM iterations \eqref{Eq:vupdate} and \eqref{Eq:supdate}.

This decentralized SDP approach has been successfully employed for monitoring large-scale power networks \cite{gg_spmag13}. To estimate the complex voltage phasor all nodes (a.k.a. power system state), measurements  are collected on real/reactive power and voltage magnitude, all of which have quadratic dependence on the unknown states. Gauss-Newton iterations have been the `workhorse' tool
for this nonlinear estimation problem; see e.g., \cite{SE_book,Wollenberg-book}. However, the iterative linearization therein could suffer from convergence issues and local optimality, especially due to the increasing variability in power grids with high penetration of renewables. With improved communication capabilities, decentralized state estimation among multiple control centers has attracted growing interest; see Fig. \ref{fig:D-SDP} illustrating three interconnected areas aiming to achieve the centralized estimation collaboratively.

\begin{figure}[!htb]
    \centering
    \begin{minipage}{0.39\textwidth}
         \includegraphics[width=\linewidth,height=3cm]{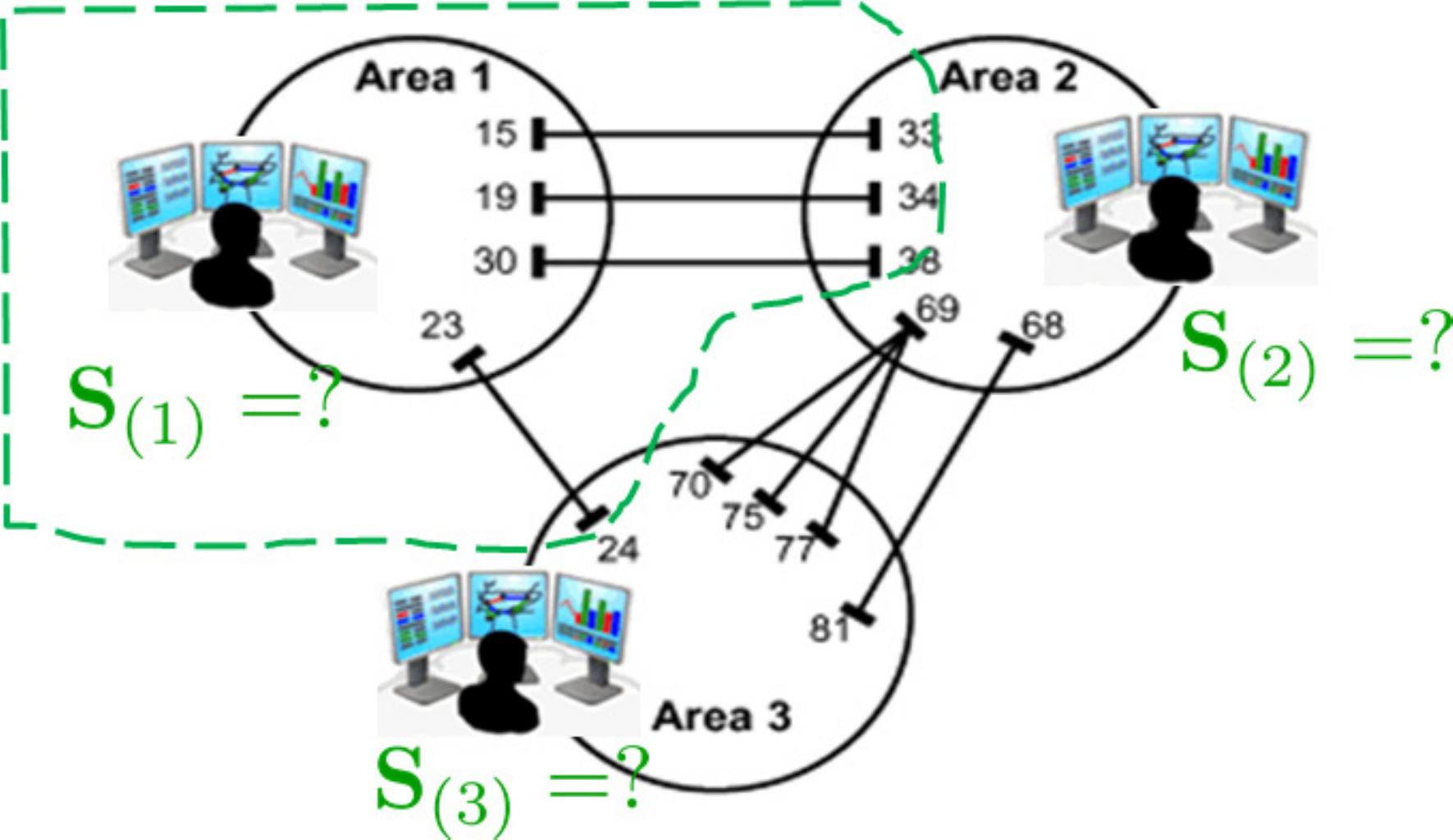}
    \end{minipage}
    \begin{minipage}{0.59\textwidth}
         \includegraphics[width=\linewidth,height=4cm]{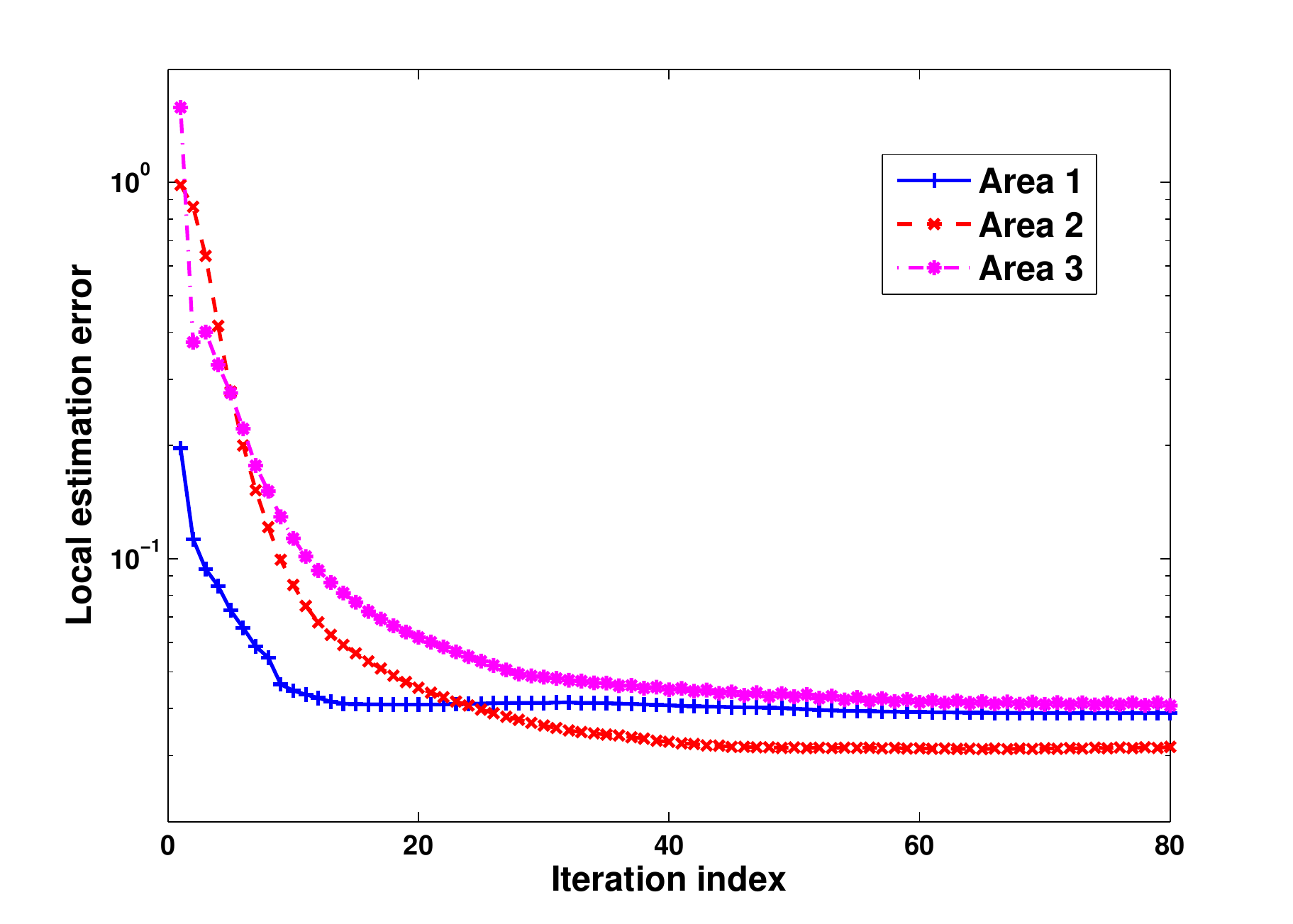}
    \end{minipage}
    \caption{(Left:) Schematic of collaborative power system state estimation among control centers of three interconnected networks (IEEE 118-bus test case). (Right:) Local state estimation error vs. iteration number using the decentralized SDP-based state estimation method. }
    \label{fig:D-SDP}
\end{figure}

A decentralized SDP-based state estimator has been developed in \cite{hzgg_jstsp14} with
reduced complexity compared to \eqref{Eq:cen_SDP}. The resultant algorithm involves only internal voltages and those of next-hop neighbors in the local matrix $\bS_{(i)}$; e.g., in Fig. \ref{fig:D-SDP} $\bS_{(1)}$ is identified by the dashed lines. Interestingly, the positive-semidefiniteness constraint for the overall $\bS$ decouples nicely into that of all local $\{\bS_{i}\}$, and the estimation error converges to the centralized performance within only a dozen iterations. The decentralized SDP framework has successfully addressed a variety of power system operational challenges, including a distributed microgrid optimal power flow solver in \cite{edhzgg_tsg13}; see also \cite{gg_spmag13} for a tutorial overview of these applications.

\subsection{Decentralized Inference}
\label{subsec:22}

Along with decentralized signal parameter estimation, a variety of inference
tasks become possible by relying on the collaborative sensing and computations
performed by networked nodes. In the special context of resource-constrained WSNs
deployed to determine the common messages broadcast by a wireless AP, the
relatively limited node reception capability makes it desirable to design a
decentralized \textit{detection} scheme for all sensors to attain sufficient
statistics for the \textit{global} problem. Another exciting application of
WSNs is environmental monitoring for e.g., inferring the presence or
absence of a pollutant over a geographical area. Limited by the local sensing
capability, it is important to develop a decentralized \textit{learning}
framework such that all sensors can \textit{collaboratively} approach the
performance as if the network wide data had been available everywhere (or at a FC
for that matter). Given the diverse inference tasks, the challenge becomes how
to design the best inter-node information exchange schemes that would allow for
minimal communication and computation overhead in specific applications.

\subsubsection{Decentralized Detection}

\runinhead{Message decoding.} A decentralized detection framework is introduced
here for the message decoding task,
which is relevant for diverse wireless communications and networking scenarios.
Consider an AP broadcasting a $p\times 1$ coded block $\bs$ to a network of
sensors, all of which know the codebook $\ccalC$ that $\bs$ belongs to. For
simplicity assume \emph{binary} codewords, and that each node $i=1,\ldots,n$
receives a same-length block of symbols $\by_i$ through a discrete, memoryless,
symmetric channel that is conditionally independent across sensors. Sensor $i$
knows its local channel from the AP, as characterized by the conditional pdf
$p(y_{il}|s_l)$ per bit $l$. Due to conceivably low signal-to-noise-ratio (SNR)
conditions, each low-cost sensor may be unable to reliably decode the message.
Accordingly, the need arises for information exchanges among single-hop neighboring
sensors to achieve the global (that is, centralized) error performance.
Given $\by_i$ per sensor $i$, the assumption on memoryless and independent
channels yields the centralized \textit{maximum-likelihood} (ML)  decoder as
\begin{equation}
    \hat{\bs}^{DEC}= \arg \max_{\bs \in \ccalC}
    p(\{\by_i\}_{i=1}^n|\bs) = \arg \min_{\bs \in \ccalC} \textstyle \sum_{l=1}^p
\sum_{i=1}^n \left[-\log p(y_{il}|s_l)\right]. \label{cendec}
\end{equation}
ML decoding amounts to deciding the most likely codeword among multiple
candidate ones and, in this sense, it can be viewed as a test of multiple
hypotheses. In this general context, belief propagation approaches have been
developed in \cite{sas06tsp},  so that all nodes can cooperate to learn
the centralized likelihood per hypothesis. However, even for linear binary
block codes, the number of hypotheses, namely the cardinality of $\ccalC$, grows
exponentially with the codeword length. This introduces high communication and
computation burden for the low-cost sensor designs.

The key here is to extract minimal sufficient statistics for the centralized
decoding problem. For binary codes, the log-likelihood terms in \eqref{cendec}
become $\log p(y_{il}|s_l)  = - \gamma_{il} s_l + \log p(y_{il}|s_l=0)$, where
\begin{equation}\label{defgamma}
    \gamma_{il} := \log \left( \frac{p(y_{il}|s_l=0)}{p(y_{il}|s_l=1)}\right)
\end{equation}
is the local log-likelihood ratio (LLR) for the bit $s_l$ at sensor $i$.
Ignoring all constant terms $\log p(y_{il}|s_l = 0)$, the ML decoding objective
ends up only depending on the sum LLRs, as given by $\hat{\bs}_{ML}= \arg \min_{\bs
\in \ccalC} \sum_{l=1}^p (\sum_{i=1}^n \gamma_{il}) s_l$. Clearly, the
sufficient statistic for solving \eqref{cendec} is the sum of all local
LLR terms, or equivalently, the average $\bbargamma_{l} = (1/n) \sum_{i=1}^n
\gamma_{il}$ for each bit $l$. Interestingly, the average of
$\{\gamma_{il}\}_{i=1}^n$ is one instance of the BLUE discussed in
Section \ref{subsubsec:DBLUE} when $\bm{\Sigma}_{y,i}=\bH_j=\bI_{p\times p}$,
since
\begin{equation}
    \bbargamma_{l}  = \textstyle\arg\min_{\gamma} \sum_{i=1}^n (\gamma_{il}
-\gamma)^2.
\end{equation}
This way, the ADMM-based decentralized learning framework in Section
\ref{sec:ADMM} allows for all sensors to collaboratively attain the sufficient
statistic for the decoding problem \eqref{cendec} via in-network processing.
Each sensor only needs to estimate a vector of the codeword length $p$, which
bypasses the exponential complexity under the framework of belief propagation.
As shown in \cite{hzggac08tsp}, decentralized \textit{soft} decoding is also
feasible since the {\it a posteriori probability (APP)} evaluator also relies on
LLR averages which are sufficient statistics, where extensions to non-binary
alphabet codeword constraints and random failing inter-sensor links are also
considered.

\begin{figure}[tb]
%\sidecaption[t]    %
\centering
\includegraphics[width=7.5cm,height=4.5cm]{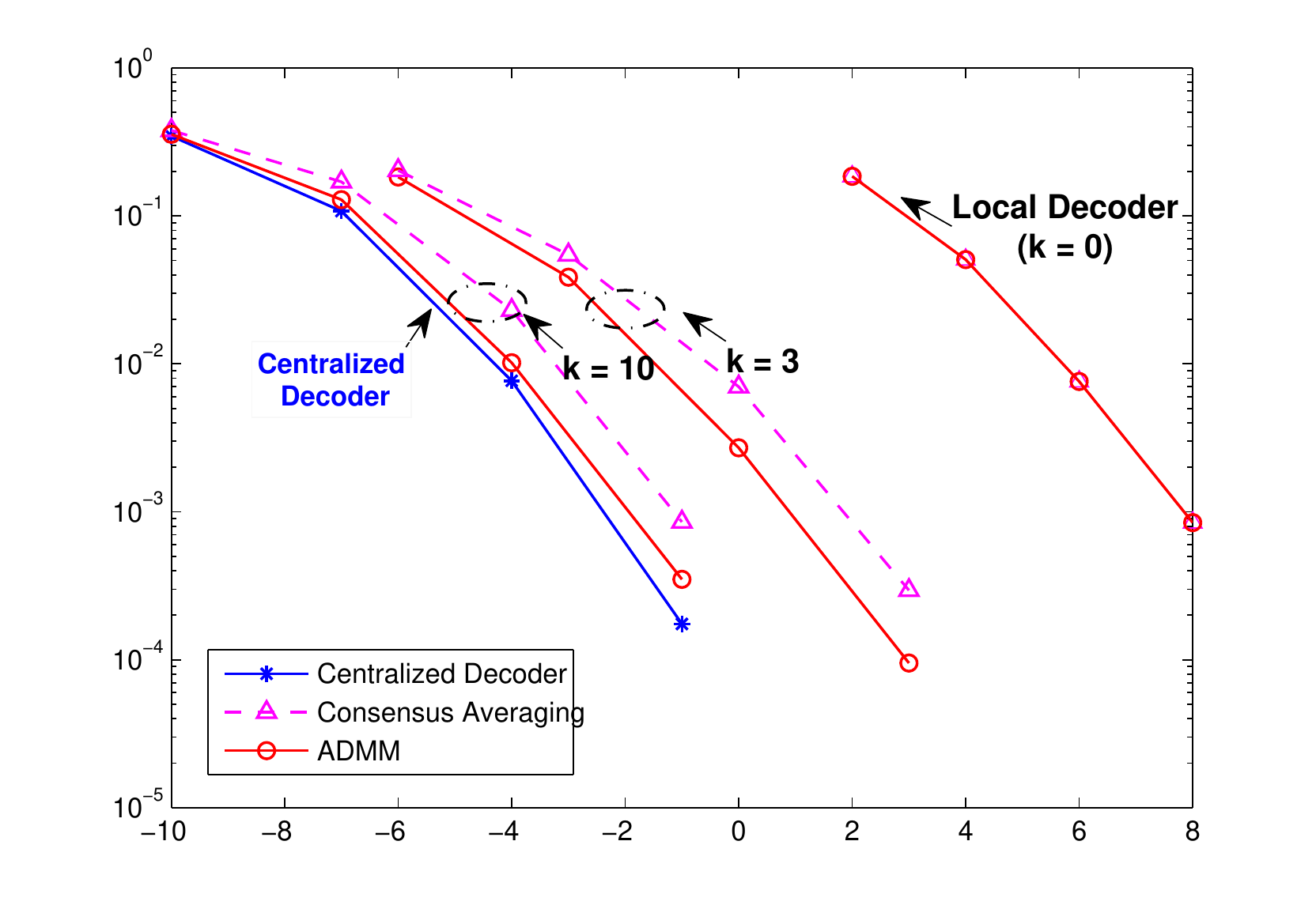}
\caption{BER vs. SNR (in dB) curves depicting the local ML decoder   vs. the
consensus-averaging decoder vs. the ADMM-based
        approach vs. the centralized ML decoder benchmark.}
    \label{fig:DiVA}
\end{figure}

The bit error rate (BER) versus SNR plot in Fig.~\ref{fig:DiVA} demonstrates the
performance
of ADMM-based in-network decoding of a convolutional code with $p=60$ and
$|\ccalC|=40$.
This numerical test involves $n=10$ sensors and AWGN AP-sensor channels
with $\sigma_i^2=10^{-SNR_i/10}$. Four schemes are compared: (i) the local ML
decoder
based on per-sensor data only (corresponds to the curve marked as $k=0$ since it
is
used to initialize the decentralized iterations); (ii) the centralized benchmark
ML
decoder (corresponds to $k=\infty$); (iii) the in-network decoder which forms
$\bbargamma_{l}$
using ``consensus-averaging'' linear iterations~\cite{Consensus_Averaging};
and, (iv) the ADMM-based decentralized algorithm. Indeed, the ADMM-based decoder
exhibits
faster convergence than its consensus-averaging counterpart; and surprisingly,
only 10 iterations suffice to
bring the decentralized BER very close to the centralized performance.

\runinhead{Message demodulation.} In a related detection scenario the common AP
message $\bs$ can be mapped to a space-time matrix, with each entry drawn from
a finite alphabet $\ccalA$. The received
block $\by_i$ per sensor $i$ typically admits a linear input/output relationship
$\by_i = \bbH_i ~\bs + \bbepsilon_i$.
Matrix $\bbH_i$ is formed from the fading AP-sensor channel, and $\bbepsilon_i$
stands for the additive white Gaussian noise of unit variance, that is assumed
uncorrelated across sensors. Since low-cost sensors have very limited budget on
number of antennas compared to the AP, the length of $\by_i$ is much shorter
than $\bs$ (i.e., $m_i<p$). Hence, the local linear demodulator using $\{\by_i,
\bbH_i\}$ may not even be able to identify $\bs$. Again, it is critical for each
sensor $i$ to cooperate with its neighbors to collectively form the global ML
demodulator
\begin{equation}
    \hhatbbs^{DEM} = \arg \max_{\bs\in\ccalA^{N}} \textstyle
-\!\sum_{i=1}^n\|\by_i-\bbH_i\bs\|^2 \!\displaystyle
    =\arg \max_{\bs\in\ccalA^{N}}\!\textstyle\left\{\!2
    \!\Big(\sum_{i=1}^n\br_i \Big)^{\top}\!\!\!\bs\! -\!
    \bs^\top\!\Big(\sum_{i=1}^n \bR_i\Big)\!\bs\!\right\}
    \label{cenML}
\end{equation}
where $\br_i := \bbH_i^\top\by_i$ and $\bR_i:= \bbH_i^\top\bbH_i$ are the sample
(cross-)covariance terms. To solve \eqref{cenML} locally, it suffices for each
sensor to acquire the network-wide average of $\{\br_i\}_{i=1}^n$, as well as
that of $\{\bR_i\}_{i=1}^n$, as both averages constitute the minimal sufficient
statistics for the centralized demodulator. Arguments similar to decentralized
decoding lead to ADMM iterations that (as with BLUE) attain locally these average
terms. These iterations constitute a viable decentralized demodulation method, whose
performance analysis in~\cite{hzacgg10twc} reveals that its error diversity
order can approach the centralized one within only a dozen of iterations.

As demonstrated by the decoding and demodulation tasks, the cornerstone of
developing a decentralized detection scheme is to extract the minimal sufficient
statistics for the centralized hypothesis testing problem. This leads to
significant complexity reduction in terms of communications and computational
overhead.

\subsubsection{Decentralized Support Vector Machines}

The merits of support vector machines (SVMs) in a centralized setting have been
well documented in various supervised classification tasks including
surveillance, monitoring, and segmentation, see e.g.,~\cite{smola}. These
applications often call for \emph{decentralized supervised learning} solutions,
when limited training data are acquired at different locations and a central
processing unit is costly or even
discouraged due to, e.g., scalability, communication overhead, or
privacy reasons. Noteworthy examples include WSNs for environmental or
structural health monitoring, as well as diagnosis of medical conditions from
patient's records distributed at different hospitals.

In this in-network classification task, a labeled
training set $\ccalT_i:=\{(\bbx_{il},y_{il})\}$ of size $m_i$ is available per
node $i$, where $\bbx_{il}\in\mathbb{R}^p$ is the input data vector
and $y_{il}\in\{-1,1\}$ denotes its corresponding class
label. Given all network-wide training data $\{\ccalT_i\}_{i=1}^n$, the
\emph{centralized} SVM seeks a maximum-margin linear discriminant function
$\hat{g}(\bbx)=\bbx^\top\hat{\bbs}+\hat{b}$, by solving the following convex
optimization problem~\cite{smola}
\begin{equation}
    \begin{aligned}
        \{\hat{\bbs},\hat{b}\}=\displaystyle\arg\min_{\bbs,~ b, \{\xi_{il}\}}\quad &
         \frac{1}{2}\left\|\bbs\right\|^2+C\sum_{i=1}^n\sum_{l=1}^{m_i} \xi_{il}& &\\
        \textrm{s. to}\quad & y_{il}(\bbs^\top\bbx_{il}+b)\geq 1-\xi_{il},&
i=1,\ldots,n&,\;l=1,\ldots,m_i\\
        &\xi_{il}\geq 0,&  i=1,\ldots,n&,\;l=1,\ldots,m_i
    \end{aligned}
    \label{SVM_P1}
\end{equation}
where the slack variables $\xi_{il}$ account for non-linearly separable training
sets, and $C$ is a tunable positive scalar that allows for controlling model
complexity. Nonlinear discriminant
functions $g(\bbx)$ can also be accommodated after mapping input vectors
$\bbx_{il}$ to a higher- (possibly infinite)-dimensional space using
e.g., kernel functions, and pursuing a generalized maximum-margin linear
classifier as in \eqref{SVM_P1}. Since the SVM classifier \eqref{SVM_P1} couples the
local datasets, early \textit{distributed} designs either rely on a centralized
processor so they are not decentralized~\cite{van08dpsvm}, or, their performance
is not guaranteed to reach that of the centralized SVM~\cite{navia06dsvm}.

A fresh view of decentralized SVM classification is taken
in~\cite{pfacgg10jmlr}, which reformulates \eqref{SVM_P1} to estimate the
parameter pair $\{\bbs,b\}$ from all local data $\ccalT_i$ after eliminating
slack variables $\xi_{il}$, namely
\begin{equation}
    \{\hat{\bbs},\hat{b}\}=\arg\min_{\bbs,~ b}~
    \frac{1}{2}\left\|\bbs\right\|^2
    +C\sum_{i=1}^n\sum_{l=1}^{m_i} \max\{0,1-y_{il}(\bbs^\top\bbx_{il}+b)\}.
    \label{SVM_P2}
\end{equation}
Notice that \eqref{SVM_P2} has the same decomposable structure that the general
decentralized learning task in \eqref{Eq:Centr_Est},
upon identifying the local cost
$f_i(\bar{\bbs};\bby_i)=\frac{1}{2n}\left\|\bbs\right\|^2
+C\sum_{l=1}^{m_i} \max\{0,1-y_{il}(\bbs^\top\bbx_{il}+b)\}$, where
$\bar{\bbs}:=[\bbs^\top,b^\top]^\top$, and
$\bby_i:=[y_{i1},\ldots,y_{im_i}]^\top$. Accordingly, all network nodes can
solve \eqref{SVM_P2} in a decentralized fashion via iterations obtained
following the ADMM-based algorithmic framework of Section \ref{sec:ADMM}. Such
a decentralized ADMM-DSVM scheme is provably convergent to the centralized SVM
classifier \eqref{SVM_P1}, and can also incorporate nonlinear discriminant
functions as detailed in~\cite{pfacgg10jmlr}.

\begin{figure}[tb]
%   \sidecaption[t]
\centering
    \includegraphics[width=7.5cm,height=4.5cm]{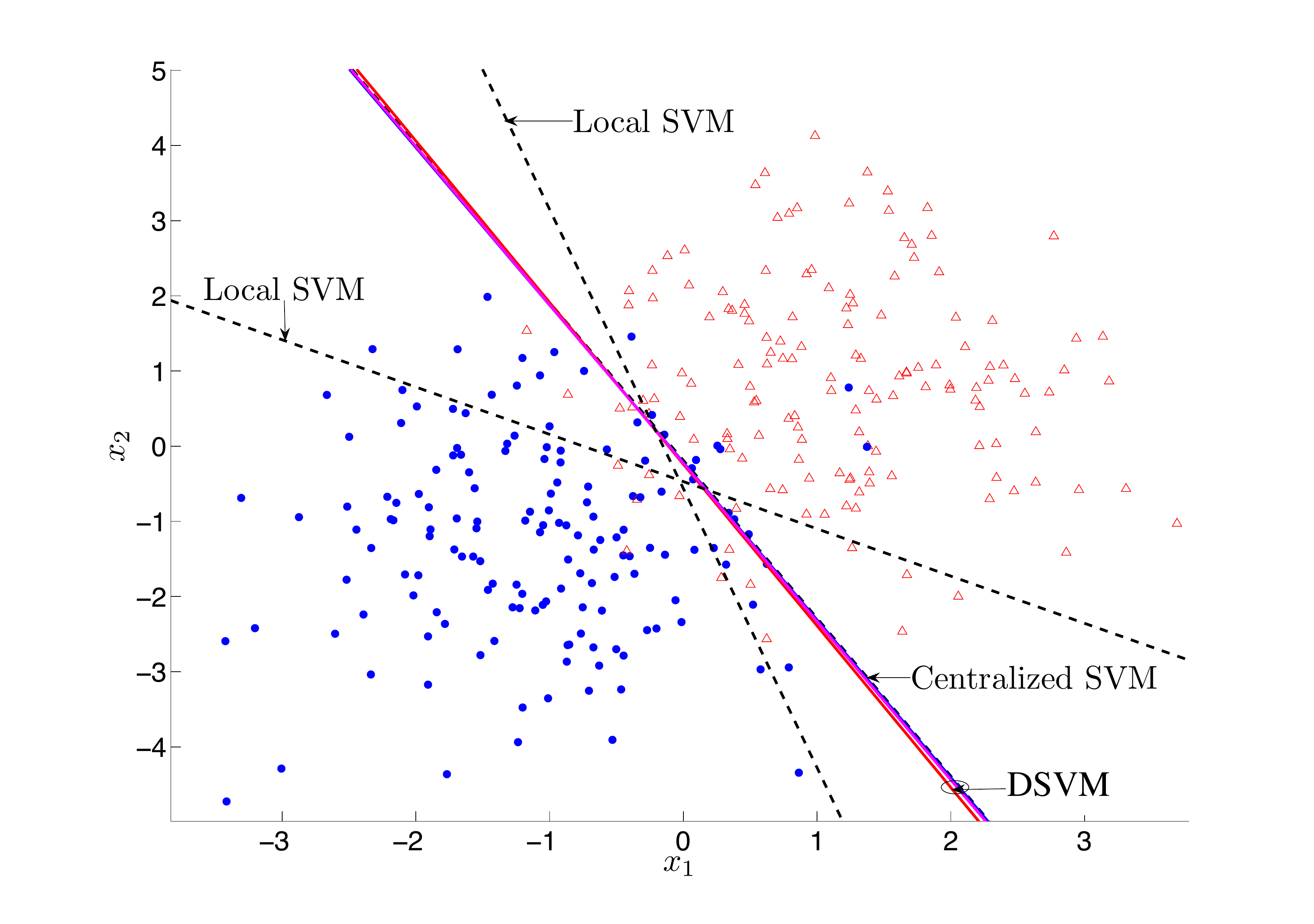}
    \caption{Decision boundary comparison among ADMM-DSVM, centralized SVM and
local SVM results for synthetic data generated from two Gaussian classes, and a
network of $n=30$ nodes.}
    \label{fig:D-SVM}
\end{figure}

To illustrate the performance of the ADMM-DSVM algorithm in~\cite{pfacgg10jmlr},
consider a randomly generated
network with $n=30$ nodes. Each node acquires labeled training examples from
two different classes, which are equiprobable and consist of random vectors
drawn from a two-dimensional (i.e., $p=2$) Gaussian distribution with common
covariance matrix $\bm{\Sigma}_x=[1, \;0;\;0, \;2]$, and mean
vectors $\bbmu_1=[-1,\; -1]^\top$ and $\bbmu_2=[1,\;
1]^\top$, respectively. The Bayes optimal classifier for this 2-class problem
is linear~\cite[Ch. 2]{duda}. To visualize
this test case, Fig.~\ref{fig:D-SVM} depicts the global training
set, along with the linear discriminant functions found by the
centralized SVM \eqref{SVM_P1} and the ADMM-DSVM at two different nodes after
400 iterations. Local SVM results for two different nodes are also included for
comparison. It is apparent
that ADMM-DSVM approaches the decision rule of its centralized counterpart,
whereas local classifiers
deviate since they neglect most of the training examples in the network.

\subsubsection{Decentralized Clustering}

Unsupervised learning using a network of wireless sensors as an exploratory
infrastructure is well motivated for inferring hidden structures in distributed data collected
by
the sensors. Different from supervised SVM-based classification tasks, each node
$i=1,\ldots,n$ has
available a set of \emph{unlabeled} observations
$\ccalX_i:=\{\bbx_{il},\;l=1,\ldots,m_i\}$, drawn from a total of $K$ classes.
In this network setting, the goal is to design local clustering rules assigning
each $\bbx_{il}$ to a cluster $k\in\{1,\ldots,K\}$. Again,
the desiderata is a decentralized algorithm capable of attaining the performance
of a benchmark clustering scheme, where all $\{\ccalX_i\}_{i=1}^n$ are centrally available for joint processing.

Various criteria are available to quantify similarity among observations in a
centralized setting, and a popular selection is the deterministic partitional
clustering (DPC) one entailing prototypical elements  (a.k.a. cluster centroids)
per class in order to avoid comparisons between every pair of observations.  Let
$\bbmu_k$ denote the prototype element
for class $k$, and $\nu_{ilk}$ the membership coefficient
of $\bbx_{il}$ to class $k$. A natural clustering problem amounts to
specifying the family of $K$ clusters with centroids $\{\bbmu_k\}_{k=1}^K$, such
that the sum of squared-errors is minimized; that is
\begin{equation}
     \min_{\{\nu_{ilk}\in\ccalV\},\{\bbmu_k\}}\sum_{i=1}^n\sum_{l=1}^{m_i}
    \sum_{k=1}^K \nu_{ilk}^\rho\left\|\bbx_{il}-\bbmu_k \right\|^2
    \label{hcluster}
\end{equation}
where $\rho \geq 1$ is a tuning parameter, and  $\ccalV:=\{\nu_{ilk}:\sum_{k}
\nu_{ilk}^\rho=1,
\;\nu_{ilk}\in[0,1],\;\forall i,l \}$ denotes the convex set of constraints on
all membership coefficients.
With $\rho=1$ and  $\{\bbmu_k\}$ fixed, \eqref{hcluster} becomes a linear
program in
$\nu_{ilk}$. Consequently, \eqref{hcluster} admits binary $\{0,1\}$ optimal
solutions giving rise to
the so-termed {\it hard} assignments, by choosing the cluster $k$ for
$\bbx_{il}$ whenever $\nu_{ilk}=1$. Otherwise,
for $\rho>1$ the optimal coefficients generally result in {\it soft} membership
assignments, and the optimal cluster is $k^* := \arg\max_{k}\nu_{ilk}^\rho$ for
$\bbx_{il}$. In either case, the DPC clustering problem \eqref{hcluster} is
NP-hard, which motivates the (suboptimal) K-means algorithm that, on a per
iteration basis, proceeds in two-steps to minimize the cost in \eqref{hcluster}
w.r.t.:
(S1) $\ccalV$ with $\{\bbmu_k\}$ fixed; and (S2) $\{\bbmu_k\}$ with $\ccalV$
fixed~\cite{lloyd82PCM}.
Convergence of this two-step alternating-minimization scheme is guaranteed
at least to a local minimum.
Nonetheless, K-means requires central availability of global information (those
variables that are fixed per step), which challenges
in-network implementations. For this reason, most early attempts are either
confined to
specific communication network topologies, or, they offer no closed-form local
solutions; see e.g.,~\cite{Nowak03dem,whk08ICML}.

To address these limitations,~\cite{pfacgg11jstsp} casts
\eqref{hcluster} [yet another instance of \eqref{Eq:Centr_Est}] as a
decentralized estimation problem. It is thus possible to leverage ADMM
iterations and solve \eqref{hcluster} in a decentralized fashion through
information exchanges among single-hop neighbors only. Albeit the non-convexity
of \eqref{hcluster}, the decentralized DPC iterations in~\cite{pfacgg11jstsp}
provably approach a local minimum arbitrarily closely, where the asymptotic
convergence holds for hard K-means with $\rho = 1$. Further extensions
in~\cite{pfacgg11jstsp} include a decentralized expectation-maximization
algorithm for probabilistic partitional clustering, and methods to handle
unknown number of classes.

\begin{figure}[t]
    \centering
    \begin{minipage}{0.49\textwidth}
        \includegraphics[width=\linewidth]{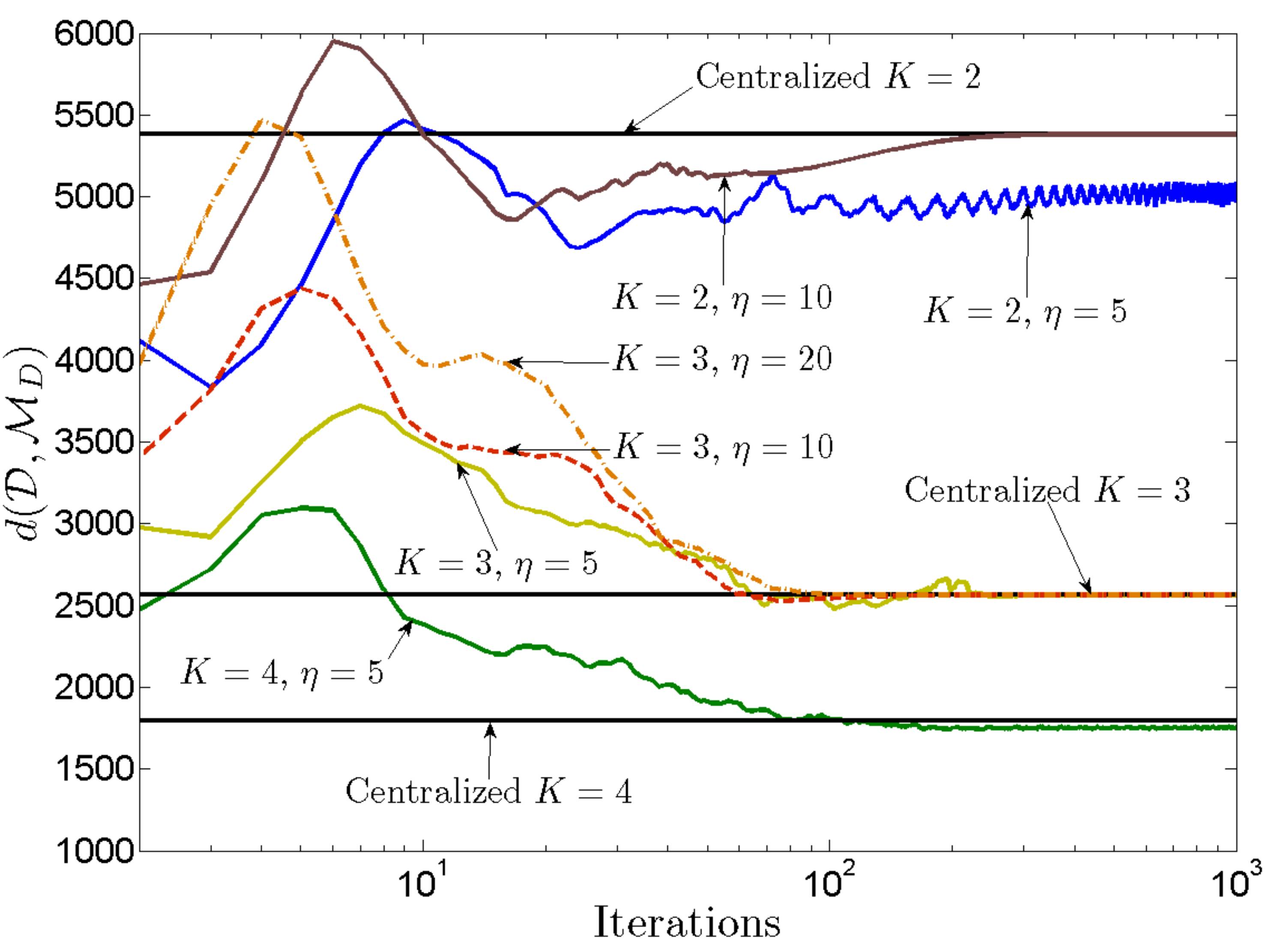}
    \end{minipage}
    \begin{minipage}{0.49\textwidth}
        \includegraphics[width=\linewidth]{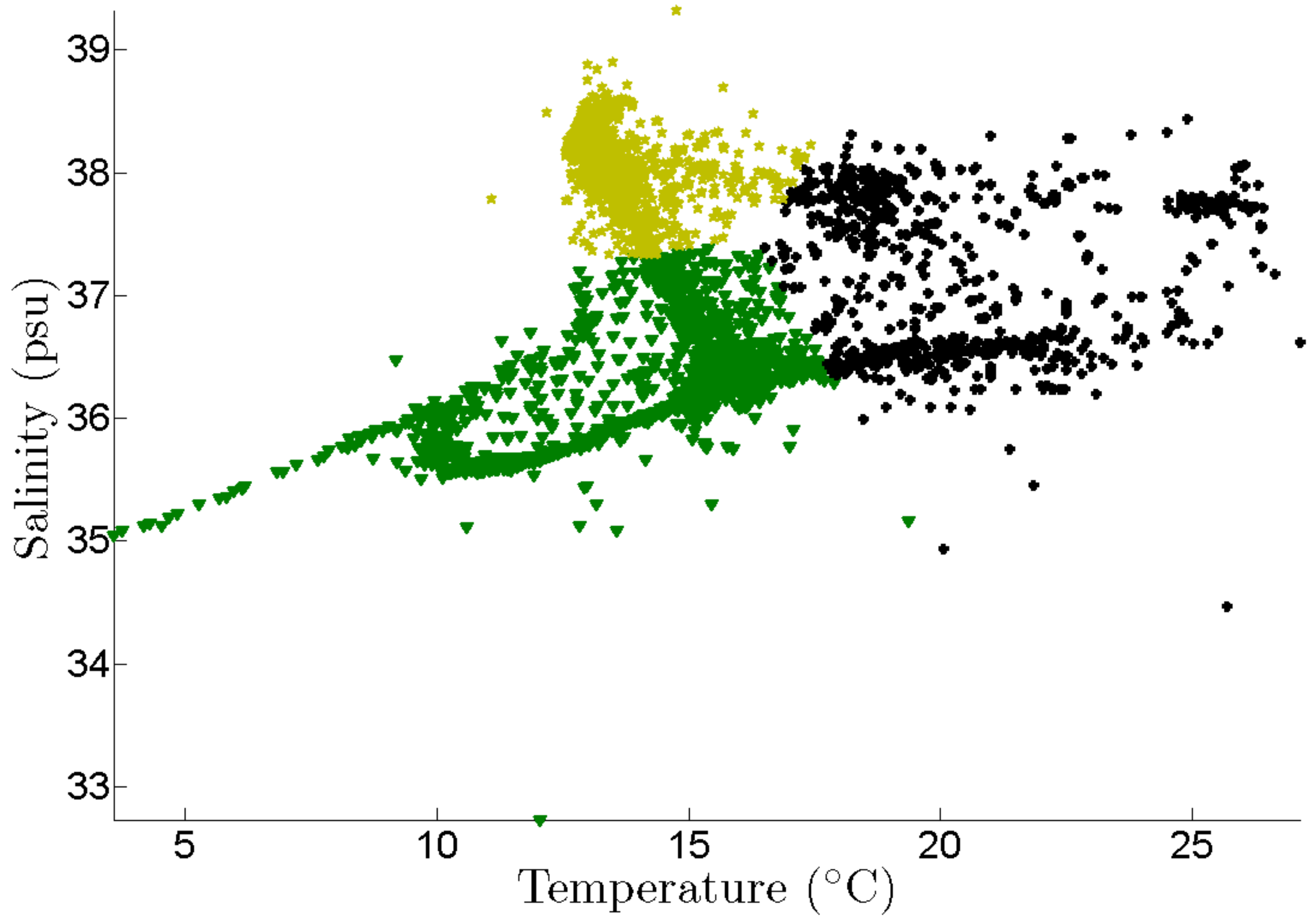}
    \end{minipage}
    \caption{Average performance of hard-DKM on a real data set using
a
        WSN with $n=20$ nodes for various values of $\eta$ and $K$ (left).
        Clustering with $K=3$ and $\eta=5$ (right) at $k=400$ iterations.}
    \label{fig:D-KM}
\end{figure}

\runinhead{Clustering of oceanographic data.} Environmental monitoring is a
typical application of WSNs. In WSNs deployed for oceanographic
monitoring, the cost of computation per node is lower than
the cost of accessing each node's observations~\cite{oceansensors}.
This makes the option of centralized processing less attractive,
thus motivating decentralized processing. Here we test the
decentralized DPC schemes of~\cite{pfacgg11jstsp} on real data
collected by multiple
underwater sensors in the Mediterranean coast of Spain~\cite{WOD}, with the goal
of
identifying regions sharing common physical characteristics. A total
of $5,720$ feature vectors were selected, each having entries the
temperature ($^\circ$C) and salinity (psu) levels ($p=2$). The
measurements were normalized to have zero mean, unit variance,
and they were grouped in $n=20$ blocks (one per sensor) of $m_i=286$
measurements each. The algebraic connectivity of the WSN is 0.2289 and the
average degree per node is 4.9. Fig. \ref{fig:D-KM} (left) shows the
performance of 25 Monte Carlo runs for the hard-DKM algorithm with
different values of the parameter $c:=\eta$. The best average
convergence rate was
obtained for $\eta=5$, attaining the average
centralized performance after 300 iterations. Tests with different
values of $K$ and $\eta$ are also included in Fig. \ref{fig:D-KM}
(left) for comparison. Note that for $K=2$ and $\eta=5$ hard-DKM
hovers around a point without converging. Choosing a larger
$\eta$ guarantees convergence of the algorithm to a unique solution.
The clustering results of hard-DKM at $k=400$ iterations for
$\eta=5$ and $K=3$ are depicted in Fig. \ref{fig:D-KM} (right).

\section{Decentralized Adaptive Estimation}\label{sec:Adaptive}
Sections \ref{sec:ADMM} and \ref{sec:Batch} dealt with  decentralized
\emph{batch} estimation, whereby
network nodes acquire data only once and then locally exchange messages to reach
consensus on
the desired estimators. In many applications however, networks are deployed to
perform estimation in a
constantly changing environment without having available a
complete statistical description of the underlying processes of
interest, e.g., with time-varying thermal or seismic sources. This motivates the
development of decentralized adaptive
estimation schemes, where nodes collect data sequentially in time and local
estimates are recursively refined ``on-the-fly.''
In settings where statistical state models are available, it is
prudent to develop
model-based tracking approaches implementing in-network Kalman or particle
filters. Next,  Section \ref{sec:ADMM}'s scope is broadened  to
facilitate real-time (adaptive) processing of network data, when the local costs
in \eqref{Eq:Centr_Est} and unknown parameters are allowed to vary
with time.

\subsection{Decentralized Least-Mean Squares}\label{sec:D_LMS}
A decentralized least-mean squares (LMS) algorithm is developed here for
adaptive estimation of (possibly) nonstationary parameters, even when
statistical information such as ensemble data covariances
are unknown. Suppose network nodes are deployed to estimate a signal vector
$\mathbf{s}(t)\in\mathbb{R}^{p\times 1}$ in a collaborative fashion subject to
single-hop communication constraints, by resorting to the linear LMS
criterion, see e.g., \cite{sk95book,Diffusion_LMS,SMG_D_LMS}.
Per time instant
$t=0,1,2,\ldots$, each node has available a regression vector
$\mathbf{h}_{i}(t)\in\mathbb{R}^{p\times 1}$ and acquires a scalar
observation $y_{i}(t)$, both assumed zero-mean without loss of
generality. Introducing the global vector
$\mathbf{y}(t):=\left[\mathbf{y}_{1}(t)\ldots
\mathbf{y}_{n}(t)\right]^\top\in\mathbb{R}^{n\times 1}$ and matrix
$\mathbf{H}(t):=\left[\mathbf{h}_{1}(t)\ldots\mathbf{h}_n(t)\right]^\top
\in\mathbb{R}^{n\times p}$, the global time-dependent LMS estimator of interest
can
be written as~\cite[p. 14]{sk95book,Diffusion_LMS,SMG_D_LMS}
%%%
\begin{equation}\label{Eq:estprblm}
    \hat{\mathbf{s}}(t):=
\arg\min_{\mathbf{s}}\mathbb{E}\left[\|\bby(t)-\mathbf{H}(t)\mathbf{s}\|^{2}\right]
={\arg}\min_{\mathbf{s}}\textstyle\sum_{i=1}^{n}\mathbb{E}\left[(y_i(t)-\mathbf{h}^\top_i(t)\mathbf{s})^{2}\right].
\end{equation}
%%%
For jointly wide-sense stationary $\{\mathbf{x}(t),\mathbf{H}(t)\}$, solving
(\ref{Eq:estprblm}) leads to the well-known Wiener filter estimate
$\hat{\mathbf{s}}_{W}=\boldsymbol{\Sigma}_{H}^{-1}\boldsymbol{\Sigma}_{Hy}$,
where
$\boldsymbol{\Sigma}_{H}:=\mathbb{E}[\mathbf{H}^\top(t)\mathbf{H}(t)]$ and
$\boldsymbol{\Sigma}_{Hy}:=\mathbb{E}[\mathbf{H}^\top(t)\mathbf{y}(t)]$; see
e.g., \cite[p.
15]{sk95book}.

For the cases where the auto- and cross-covariance matrices $\boldsymbol{\Sigma}_{H}$ and $\boldsymbol{\Sigma}_{Hy}$ are unknown, the approach
followed  here to develop
the decentralized (D-) LMS algorithm includes two main building
blocks: (i) recast \eqref{Eq:estprblm} into an equivalent form amenable
to in-network processing via the ADMM framework of Section \ref{sec:ADMM};  and
(ii) leverage stochastic approximation iterations~\cite{kushner} to obtain an
adaptive LMS-like algorithm that can handle the
unavailability/variation of statistical information. Following those algorithmic
construction steps outlined in Section \ref{sec:ADMM}, the following updating
recursions are obtained for the multipliers $\mathbf{v}_i(t)$ and the
local estimates $\mathbf{s}_i(t+1)$ at time instant $t+1$ and $i=1,\ldots,n$
\begin{align}
\mathbf{v}_i(t)&=\mathbf{v}_i(t-1)+c\sum_{j\in\ccalN_i}[\mathbf{s}_i(t)-\mathbf{s}_j(t)]
    \label{Eq:v_updates_LMS}\\
\mathbf{s}_i(t+1)&=\arg\min_{\mathbf{s}_i}\left\{\mathbb{E}\left[(y_{i}(t+1)-\mathbf{h}^\top_i(t+1)\mathbf{s}_{i})^{2}\right]
    +\mathbf{v}_i^\top(t)\mathbf{s}_{i}\phantom{\sum_{j\in{\cal N}_i}}
\right.\nonumber\\
    &\hspace{2.1cm}\left.+ c\sum_{j\in{\cal
N}_i}\left\|\mathbf{s}_i-\frac{\mathbf{s}_i(t)+\mathbf{s}_{j}(t)}{2}\right\|^2\right\}.
    \label{Eq:s_updates_LMS}
\end{align}
%%%
It is apparent that after differentiating \eqref{Eq:s_updates_LMS} and setting
the gradient equal to zero,  $\mathbf{s}_i(t+1)$ can be obtained
as the root of an equation of the form
\begin{equation}\label{Eq:Unknown_Equation}
\mathbb{E}[\boldsymbol{\varphi}(\mathbf{s}_i,y_i(t+1),\mathbf{h}_i(t+1))]=\mathbf{0}
\end{equation}
where $\boldsymbol{\varphi}$ corresponds to the stochastic gradient of the
cost in \eqref{Eq:s_updates_LMS}. However, the previous equation cannot be
solved since the nodes do not have available any statistical information about
the acquired data. Inspired by stochastic approximation
techniques (such as the celebrated Robbins-Monro algorithm; see e.g.,\cite[Ch.
1]{kushner}) which iteratively find the root of
\eqref{Eq:Unknown_Equation} given noisy
observations
$\{\boldsymbol{\varphi}(\mathbf{s}_i(t),y_i(t+1),\mathbf{h}_i(t+1))\}_{t=0}^{\infty}$,
one can just
drop the unknown expected value to obtain the
following D-LMS (i.e., stochastic gradient) updates
\begin{equation}\label{Eq:s_update_LMS}
\mathbf{s}_i(t+1)=\mathbf{s}_i(t)+\mu\left[\mathbf{h}_i(t+1)e_i(t+1)-\mathbf{v}_i(t)-
    c\sum_{j\in{\cal N}_i}[\mathbf{s}_i(t)-\mathbf{s}_{j}(t)]\right]
\end{equation}
where $\mu$ denotes a constant step-size, and $e_i(t+1):=2[
y_i(t+1)-\mathbf{h}_i^\top(t+1)\mathbf{s}_i(t)]$ is twice the local \emph{a
priori} error.

Recursions \eqref{Eq:v_updates_LMS} and \eqref{Eq:s_update_LMS} constitute the
D-LMS algorithm, which can
be viewed as a stochastic-gradient counterpart of D-BLUE in Section
\ref{subsubsec:DBLUE}. D-LMS is a pioneering approach for decentralized
online learning, which blends for the first time affordable (first-order)
stochastic approximation steps with parallel ADMM iterations.
The use of a constant step-size $\mu$ endows D-LMS with tracking
capabilities. This is desirable in a constantly changing environment,
within which e.g., WSNs are envisioned to operate. The D-LMS algorithm is stable
and converges even in the presence of
inter-node communication noise (see details in \cite{SMG_D_LMS,MSG_D_LMS}).
Further, closed-form expressions for the evolution and the
steady-state mean-square error (MSE), as well as selection guidelines for the
step-size $\mu$ can be found in \cite{MSG_D_LMS}.

\begin{figure}[t]
    \centering
    \begin{minipage}{0.49\textwidth}
        \includegraphics[width=\linewidth]{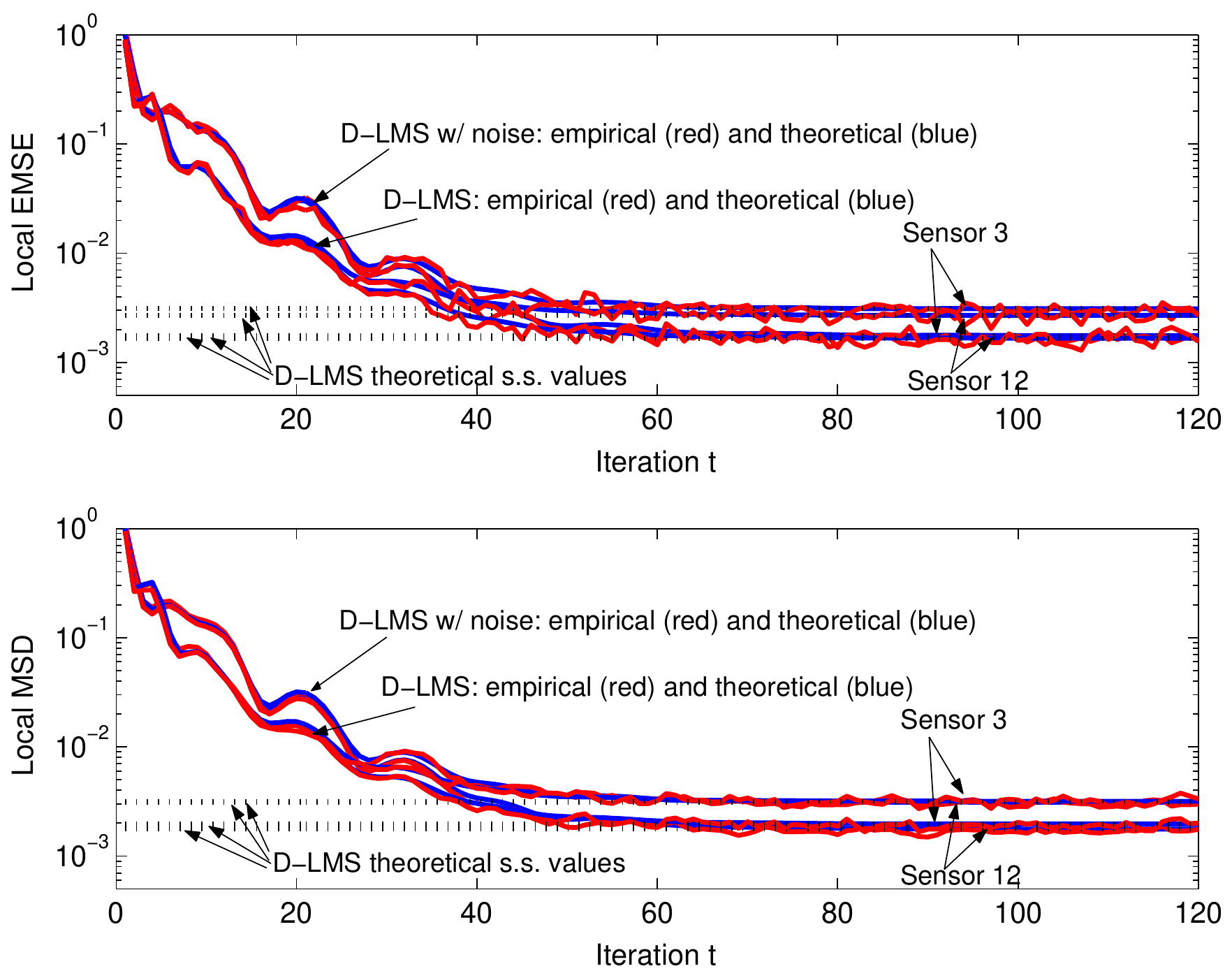}
    \end{minipage}
    \begin{minipage}{0.49\textwidth}
        \includegraphics[width=\linewidth]{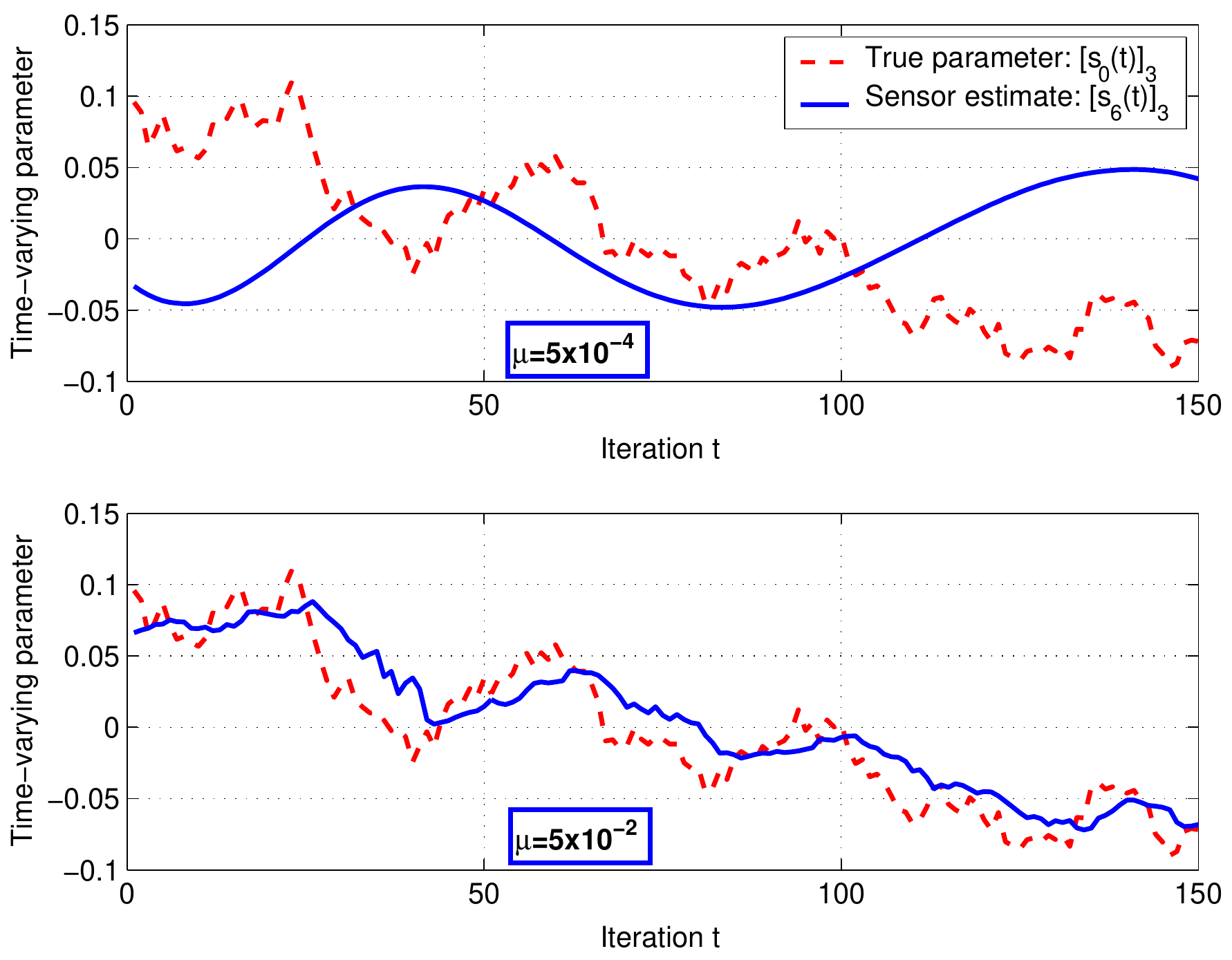}
    \end{minipage}
    \caption{Tracking with D-LMS. (left) Local MSE performance metrics both with
and without inter-node
    communication noise for sensors 3 and 12; and (right) True and estimated
time-varying parameters for a representative
    node, using slow and optimal adaptation levels.}
    \label{fig:D-LMS}
\end{figure}

Here we test the tracking performance of D-LMS with a computer simulation. For a
random geometric graph with $n=20$ nodes, network-wide observations $y_i$ are
linearly related to
a large-amplitude slowly time-varying parameter vector
$\bbs_0(t)\in\mathbb{R}^4$. Specifically, $\bbs_0(t)=\bm{\Theta}
\bbs_0(t-1)+\bm{\zeta}(t)$, where
$\bm{\Theta}=(1-10^{-4})\textrm{diag}(\theta_1,\ldots,\theta_p)$ with
$\theta_i\sim\mathcal{U}[0,1]$. The driving
noise is normally distributed with $\bm{\Sigma}_{\zeta}=10^{-4}\bbI_p$. To model
noisy links, additive white Gaussian noise
with variance $10^{-2}$ is present at the receiving end. For $\mu=5\times
10^{-2}$, Fig. \ref{fig:D-LMS} (left)
depicts the local performance of two representative nodes through the evolution
of the excess mean-square error
$\textrm{EMSE}_i(t)=\mathbb{E}[(\bbh_i^\top(t)[\bbs_i(t-1)-\bbs_0(t-1)])^2]$ and
the mean-square deviation
$\textrm{MSD}_i(t)=
\mathbb{E}[\|\bbs_i(t)-\bbs_0(t)\|^2]$ figures of merit.
Both noisy and ideal links are considered, and the empirical curves closely
follow the
theoretical trajectories derived in~\cite{MSG_D_LMS}. Steady-state limiting
values are also extremely accurate. As intuitively expected and suggested by the
analysis,
a performance penalty due to non-ideal links is also apparent.
Fig. \ref{fig:D-LMS} (right)
illustrates how the adaptation level affects the resulting per-node
estimates when tracking time-varying parameters with D-LMS. For
$\mu=5\times 10^{-4}$ (slow adaptation) and $\mu=5\times 10^{-2}$ (near optimal
adaptation),
we depict the third entry of the parameter vector $[\bbs_0(t)]_3$ and the
respective estimates
from the randomly chosen sixth node. Under optimal adaptation the local
estimate closely tracks the true variations, while -- as expected -- for
the smaller step-size D-LMS fails to provide an accurate
estimate~\cite{MSG_D_LMS,sk95book}.

\subsection{Decentralized Recursive Least-Squares}\label{sec:D_RLS}
%%%
The recursive least-squares (RLS) algorithm has well-appreciated
merits for reducing complexity and storage requirements, in online estimation of stationary signals, as well as for tracking slowly-varying nonstationary
processes~\cite{sk95book,Estimation_Theory}.
RLS is especially attractive when the state and/or data model are not
available (as with LMS), and fast convergence rates are at a premium.
Compared to the LMS scheme, RLS typically offers faster convergence and
improved estimation performance at the cost of
higher computational complexity. To enable these valuable tradeoffs in the
context of in-network processing,
the  ADMM framework of Section \ref{sec:ADMM} is utilized
here to derive a decentralized (D-) RLS adaptive scheme that can be employed for
distributed
localization and power spectrum estimation (see also~\cite{MSG_D_RLS,MG_D_RLS}
for further details
on the algorithmic construction and convergence claims).

Consider the data setting and linear regression task in Section \ref{sec:D_LMS}.
The RLS estimator
for the unknown parameter $\mathbf{s}_0(t)$ minimizes the
exponentially weighted least-squares (EWLS) cost, see e.g.,
\cite{sk95book,Estimation_Theory}
\begin{equation}\label{estprblm}
\hat{\mathbf{s}}_{\textrm{ewls}}(t):={\arg}\min_{\mathbf{s}}\sum_{\tau=0}^{t}\sum_{i=1}^{n}
     \gamma^{t-\tau}\left[y_i(\tau)-\mathbf{h}^\top_i(\tau)\mathbf{s}\right]^{2}+
    \gamma^{t}\mathbf{s}^{T}\boldsymbol{\Phi}_0\mathbf{s}
\end{equation}
where $\gamma\in(0,1]$ is a forgetting factor, while the positive
definite matrix $\boldsymbol{\Phi}_0$ is included for regularization. Note
that in forming the EWLS estimator at time $t$, the entire history of data
$\{y_i(\tau),\mathbf{h}_i(\tau)\}_{\tau=0}^{t}$ for $i=1,\ldots,n$ is
incorporated in the online estimation process. Whenever
$\gamma<1$, past data are exponentially discarded thus enabling
tracking of nonstationary processes.

Again to decompose the cost function in \eqref{estprblm}, in which
summands are coupled through the global variable $\mathbf{s}$, we
introduce auxiliary variables $\{\mathbf{s}_i\}_{i=1}^n$ that
represent local estimates per node $i$. These local
estimates are utilized to form the convex \emph{constrained} and separable
minimization problem in \eqref{Eq:Constr_Equi}, which can be solved using ADMM
to yield the following decentralized iterations (details in
\cite{MSG_D_RLS,MG_D_RLS})
\begin{align}
\mathbf{v}_i(t)&=\mathbf{v}_i(t-1)+c\sum_{j\in\ccalN_i}[\mathbf{s}_i(t)-\mathbf{s}_j(t)]\label{Eq:Vj_RLS}\\
     \mathbf{s}_i(t+1)&=\boldsymbol{\Phi}_i^{-1}(t+1)\boldsymbol{\psi}_i(t+1)
     -\frac{1}{2}\boldsymbol{\Phi}_i^{-1}(t+1)\mathbf{v}_i(t)\label{Eq:Sj_RLS}
\end{align}
where
$\boldsymbol{\Phi}_i(t+1):=\sum_{\tau=0}^{t+1}\gamma^{t+1-\tau}\mathbf{h}_i(\tau)\mathbf{h}_i^\top(\tau)+n^{-1}\gamma^{t+1}
    \boldsymbol{\Phi}_0$ and
\begin{align}\label{Eq:Phiandpsi}
     \boldsymbol{\Phi}_i^{-1}(t+1)&=\gamma^{-1}\boldsymbol{\Phi}_i^{-1}(t)-
\frac{\gamma^{-1}\boldsymbol{\Phi}_i^{-1}(t)\bbh_i(t+1)\bbh_i^\top(t+1)\boldsymbol{\Phi}_i^{-1}(t)}
     {\gamma+\bbh_i^\top(t+1)\boldsymbol{\Phi}_i^{-1}(t)\bbh_i(t+1)}\\
    \boldsymbol{\psi}_i(t+1)&
     :=\sum_{\tau=0}^{t+1}\gamma^{t+1-\tau}\mathbf{h}_i(\tau)y_i(\tau)=
    \gamma\boldsymbol{\psi}_i(t)+\mathbf{h}_i(t+1)y_i(t+1).
\end{align}
The D-RLS recursions \eqref{Eq:Vj_RLS} and \eqref{Eq:Sj_RLS} involve similar
inter-node communication exchanges
as in D-LMS. It is recommended to initialize the matrix recursion with
$\boldsymbol{\Phi}_i^{-1}(0)=n\boldsymbol{\Phi}_0^{-1}:=\delta\bbI_p$, where
$\delta>0$
is chosen sufficiently large~\cite{sk95book}. The local estimates in D-RLS
converge in the mean-sense to the true $\mathbf{s}_0$ (time-invariant
case), even when information exchanges are imperfect. Closed-form expressions
for the bounded estimation MSE along with numerical tests and comparisons with
the incremental RLS~\cite{inc_RLS} and diffusion RLS~\cite{Diffusion_RLS}
algorithms can be found in~\cite{MG_D_RLS}.

\runinhead{Decentralized spectrum sensing using WSNs.} A WSN application where
the need for linear regression arises, is
spectrum estimation for the purpose of environmental monitoring.
Suppose sensors comprising a WSN deployed over some area of interest
observe a narrowband source to determine its spectral peaks. These
peaks can reveal hidden periodicities due to e.g., a natural heat or
seismic source. The source of interest propagates through multi-path
channels and is contaminated with additive noise present at the
sensors. The unknown source-sensor channels may introduce deep fades
at the frequency band occupied by the source. Thus, having each
sensor operating on its own may lead to faulty assessments. The
available spatial diversity to effect improved spectral estimates,
can only be achieved via sensor collaboration as in the decentralized
estimation algorithms presented in this chapter.

Let ${\theta}(t)$ denote the evolution of the source signal in time, and
suppose that ${\theta}(t)$ can be modeled as an autoregressive (AR)
process~\cite[p. 106]{Stoica_Book}
\begin{equation*}
    {\theta}(t)=-\sum_{\tau=1}^{p}\alpha_{\tau}{\theta}(t-\tau)+w(t)
\end{equation*}
where $p$ is the order of the AR process, while $\{\alpha_\tau\}$
are the AR coefficients and $w(t)$ denotes driving white noise. The source
propagates to sensor $i$ via a channel modeled as an FIR filter
$C_i(z)=\sum_{l=0}^{L_i-1}c_{il}z^{-l}$, of unknown order $L_i$ and
tap coefficients $\{c_{il}\}$ and is contaminated with additive
sensing noise $\bar{\epsilon}_i(t)$ to yield the observation
\begin{equation*}
    y_i(t)=\sum_{l=0}^{L_i-1}c_{il}\theta(t-l)+\bar{\epsilon}_i(t).
\end{equation*}
Since $y_i(t)$ is an autoregressive moving average (ARMA) process,
then~\cite{Stoica_Book}
\begin{equation}\label{arma_xj_intro}
    y_i(t)=-\sum_{\tau=1}^{p}\alpha_\tau
     y_i(t-\tau)+\sum_{\tau'=1}^{m}\beta_{\tau'}\tilde{\eta}_i(t-\tau')
\end{equation}
where the MA coefficients $\{\beta_{\tau'}\}$ and the variance of the
white noise process $\tilde{\eta}_i(t)$ depend on $\{c_{il}\}$,
$\{\alpha_{\tau}\}$ and the variance of the noise terms $w(t)$ and
$\bar{\epsilon}_i(t)$. For the purpose of determining spectral peaks, the
MA term in \eqref{arma_xj_intro} can be treated as observation noise,
i.e., $\epsilon_i(t):=\sum_{\tau'=1}^{m}\beta_{\tau'}\tilde{\eta}_i(t-\tau')$.
This is very important since this way sensors do not have to know the
source-sensor
channel coefficients as well as the noise variances. Accordingly, the spectral
content
of the source can be estimated provided sensors estimate the coefficients
$\{\alpha_{\tau}\}$. To this end, let $\bbs_0:=[\alpha_1\ldots\alpha_p]^\top$ be
the
unknown parameter of interest. From \eqref{arma_xj_intro} the regression
vectors are given as $\bbh_i(t)=[-y_i(t-1)\ldots -y_i(t-p)]^\top$, and can be
acquired directly from the sensor measurements $\{y_i(t)\}$ without the need of
training/estimation.

\begin{figure}[t]
    \centering
    \begin{minipage}{0.49\textwidth}
         \includegraphics[width=\linewidth]{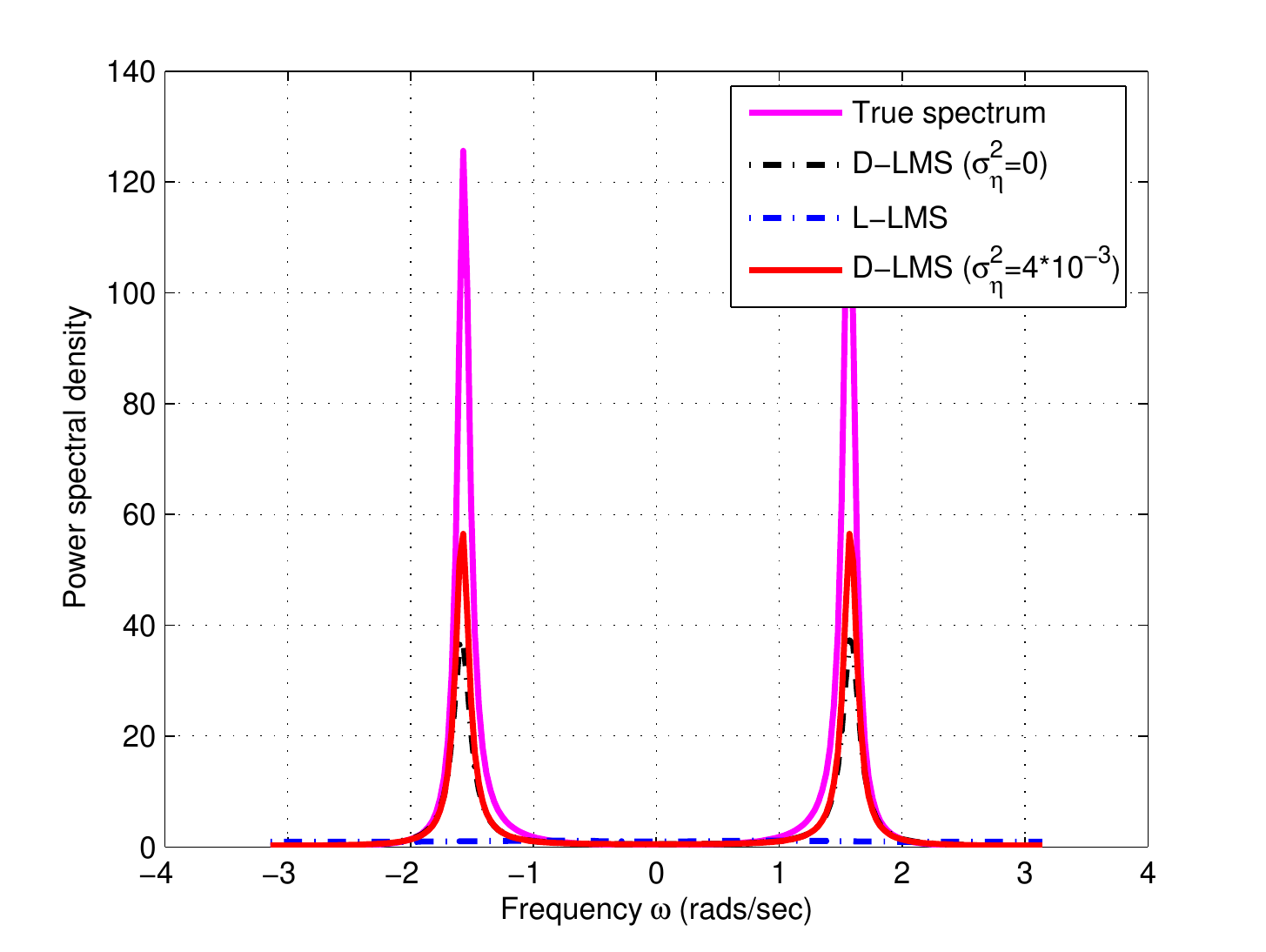}
    \end{minipage}
    \begin{minipage}{0.43\textwidth}
        \includegraphics[width=\linewidth]{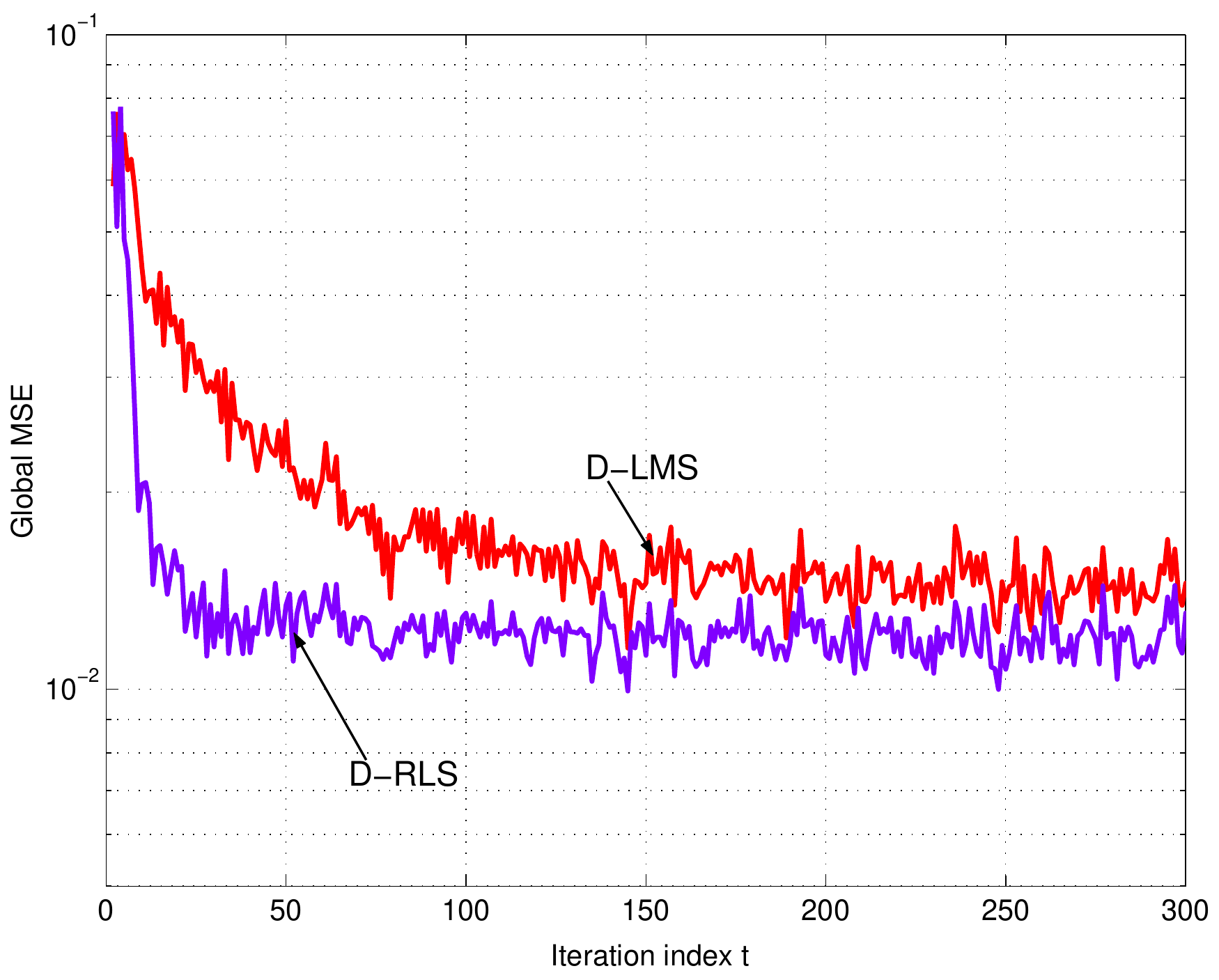}
    \end{minipage}
    \caption{D-LMS in a power spectrum estimation task. (left) The true
        narrowband spectra is compared to the estimated PSD, obtained after
        the WSN runs the D-LMS and (non-cooperative) L-LMS algorithms. The
        reconstruction results correspond to a sensor whose multipath
        channel from the source introduces a null at $\omega=\pi/2=1.57$. (right)
        Global MSE evolution (network learning curve) for the D-LMS
        and D-RLS algorithms.}
    \label{fig:spectrum_estimation_intro}
\end{figure}

Performance of the decentralized adaptive algorithms described so far
is illustrated next, when applied to the aforementioned power spectrum
estimation task. For the numerical experiments, an ad hoc WSN with
$n=80$ sensors is simulated as a realization of a random geometric
graph. The source-sensor channels corresponding to a few of the sensors are set
so that they
have a null at the frequency where the AR source has a peak, namely
at $\omega=\pi/2$. Fig. \ref{fig:spectrum_estimation_intro} (left) depicts the
actual power spectral density (PSD) of the source as well as the
estimated PSDs for one of the sensors affected by a bad channel. To
form the desired estimates in a distributed fashion, the WSN runs the
local (L-) LMS and the D-LMS algorithm outlined in Section \ref{sec:D_LMS}.
The L-LMS is a non-cooperative scheme since each sensor, say the
$i$th, independently runs an LMS adaptive filter fed by its local
data $\{y_i(t),\bbh_i(t)\}$ only. The
experiment involving D-LMS is performed under ideal and noisy inter-sensor links.
Clearly, even in the presence of communication noise D-LMS exploits
the spatial diversity available and allows all sensors to estimate
accurately the actual spectral peak, whereas L-LMS leads the
problematic sensors to misleading estimates.

For the same setup, Fig. \ref{fig:spectrum_estimation_intro} (right) shows the
global learning curve evolution
$\textrm{MSE}(t)=(1/n)\sum_{i=1}^{n}\|y_i(t)-\bbh_i^\top(t)\bbs_i(t-1)\|^2$.
The D-LMS and the D-RLS algorithms are
compared under ideal communication links. It is apparent that D-RLS
achieves improved performance both in terms of convergence rate and steady
state MSE. As discussed in Section \ref{sec:D_RLS} this comes at the price
of increased computational complexity per sensor, while the communication
costs incurred are identical.

\subsection{Decentralized Model-based Tracking}\label{sec:D_KF_KS}
The decentralized adaptive schemes in Secs. \ref{sec:D_LMS} and \ref{sec:D_RLS}
are
suitable for tracking slowly time-varying signals in settings where no
statistical models are
available. In certain cases, such as target tracking, state evolution models can
be derived
and employed by exploiting the physics of the problem. The availability of such
models
paves the way for improved state tracking via Kalman filtering/smoothing
techniques,
e.g., see \cite{Opt_Filtering_Moore,Estimation_Theory}.  Model-based
decentralized
Kalman filtering/smoothing as well as particle filtering schemes for multi-node
networks are briefly outlined here.

Initial attempts to distribute the centralized KF recursions
(see~\cite{Olfati_Kalman} and references in~\cite{sgrr08tsp}) rely on
consensus-averaging~\cite{Consensus_Averaging}. The idea is to estimate across
nodes those sufficient statistics (that are expressible
in terms of network-wide averages) required to form the corrected state and
corresponding corrected state error covariance
matrix. Clearly, there is an inherent delay in obtaining these estimates
confining the operation of such schemes only to applications with
slow-varying state vectors $\mathbf{s}_0(t)$, and/or fast communications needed
to complete multiple
consensus iterations within the time interval separating the acquisition of
consecutive measurements $y_i(t)$ and $y_i(t+1)$.
Other issues that may lead to instability in existing decentralized KF
approaches are detailed in~\cite{sgrr08tsp}.

Instead of filtering,  the delay incurred by those inner-loop consensus
iterations motivated the
consideration of fixed-lag decentralized Kalman smoothing (KS)
in~\cite{sgrr08tsp}. Matching consensus iterations
with those time instants of data acquisition, fixed-lag smoothers allow sensors
to form local MMSE optimal smoothed estimates, which take advantage of all
acquired measurements within the ``waiting period.''
The ADMM-enabled decentralized KS in~\cite{sgrr08tsp} also overcomes the
noise-related
limitations of consensus-averaging algorithms~\cite{xbk07jpdc}. In the presence
of communication noise, these estimates
converge in the mean sense, while their noise-induced variance remains bounded.
This noise resiliency allows sensors to exchange quantized data further lowering
communication cost.  For a tutorial treatment of decentralized Kalman filtering approaches using WSNs (including the decentralized ADMM-based KS
of~\cite{sgrr08tsp} and strategies to reduce the communication cost of
state estimation problems), the interested reader is referred to~\cite{dkf_control_mag}. These reduced-cost strategies exploit
the redundancy in information provided by individual observations
collected at different sensors, different observations collected at
different sensors, and different observations acquired at the same
sensor.

On a related note, a collaborative algorithm is developed
in~\cite{cg_cartography} to estimate the channel
gains of wireless links in a geographical area. Kriged Kalman filtering
(KKF)~\cite{ripley}, which is a tool with
widely appreciated merits in spatial statistics and geosciences, is adopted and
implemented in a
decentralized fashion leveraging the ADMM framework described here.
The distributed KKF algorithm requires only local message passing to track the
time-variant so-termed ``shadowing field''
using a network of radiometers, yet it provides a global view of the radio
frequency (RF) environment through consensus
iterations; see also Section \ref{ssec:rf} for further elaboration on spectrum sensing carried out via wireless cognitive radio networks.

To wrap-up the discussion, consider a network of collaborating agents (e.g.,
robots) equipped
with wireless sensors measuring distance and/or bearing
from a target that they wish to track. Even if state models are available,
the nonlinearities present in these measurements prevent sensors from employing
the
clairvoyant (linear) Kalman tracker discussed so far. In response to these
challenges,~\cite{dpf} develops a set-membership
constrained particle filter (PF) approach that: (i) exhibits performance
comparable to the centralized PF; (ii) requires only communication
of particle weights among neighboring sensors; and (iii) it can
afford both consensus-based and incremental averaging implementations.
Affordable inter-sensor communications are enabled through
a novel distributed adaptation scheme, which considerably reduces
the number of particles needed to achieve a given performance. The interested
reader is referred to~\cite{dpf_tutorial} for a recent tutorial account of
decentralized PF
in multi-agent networks.

\section{Decentralized Sparsity-regularized Rank Minimization}\label{sec:Sparse}
\label{sec:2}

Modern network data sets typically involve a large number of attributes. This
fact motivates predictive models offering a
\emph{sparse}, broadly meaning parsimonious, representation in terms of a few
attributes. Such low-dimensional models facilitate interpretability and enhanced predictive performance. In this context, this section deals with ADMM-based decentralized algorithms for sparsity-regularized rank minimization.
It is argued that such algorithms are key to unveiling Internet traffic
anomalies given ubiquitous link-load measurements.
Moreover, the notion of RF cartography is subsequently introduced to exemplify
the development of a paradigm infrastructure for situational awareness at the
physical layer of wireless cognitive radio (CR) networks. A (subsumed)
decentralized sparse linear regression algorithm
is outlined to accomplish the aforementioned cartography task.

\subsection{Network Anomaly Detection Via Sparsity and Low Rank}

Consider a backbone IP network, whose abstraction is a graph  with $n$ nodes
(routers) and  $L$ physical links. The operational goal of the network is to transport a set of $F$ origin-destination (OD) traffic flows associated with specific OD (ingress-egress router) pairs. Let $x_{l,t}$ denote
the traffic volume (in bytes or packets) passing through link $l \in
\{1,\ldots,L\}$ over a fixed time
interval $(t,t+\Delta t)$. Link counts across the entire network
are collected in the vector $\bx_t\in\mathbb{R}^L$, e.g., using the ubiquitous
SNMP protocol.
Single-path routing is adopted here, meaning a given flow's traffic is carried
through multiple
links connecting the corresponding source-destination pair along a single path.
Accordingly,
over a discrete time horizon $t \in [1,T]$ the measured link counts
$\bX:=[x_{l,t}]\in\mathbb{R}^{L \times T}$ and (unobservable) OD flow traffic
matrix
$\bZ:=[z_{f,t}]\in\mathbb{R}^{F \times T}$, are thus related through
$\bX=\bR\bZ$~\cite{lakhina}, where
the so-termed routing matrix $\bR:=[r_{l,f}]\in\{0,1\}^{L\times F}$
is such that $r_{l,f} = 1$ if link $l$ carries the flow $f$, and
zero otherwise.
The routing matrix is `wide,' as for backbone networks the number
of OD flows is much larger than the number of physical links $(F\gg L)$.
A cardinal property of the traffic matrix is noteworthy.
Common temporal patterns across OD traffic flows in addition to their almost
periodic behavior, render most rows (respectively columns) of the traffic
matrix linearly dependent, and thus $\bZ$
typically has \textit{low rank}. This intuitive property has been extensively
validated
with real network data; see Fig. \ref{fig:fig_flows} and e.g.,~\cite{lakhina}.

\begin{figure}[ht]
%\sidecaption[t]
\centering
\includegraphics[width=7.5cm]{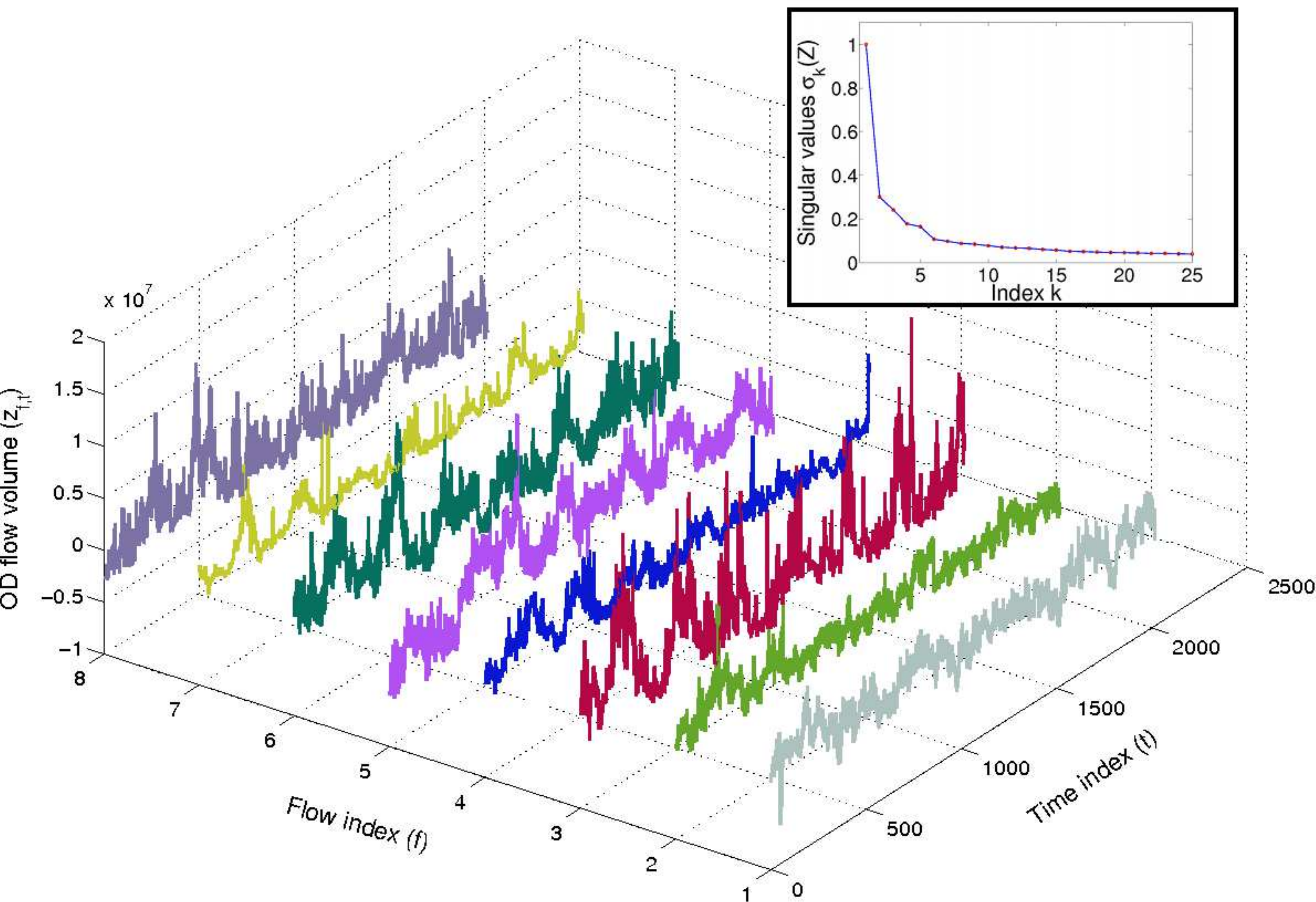}
\caption{Volumes of $6$ representative (out of $121$ total) OD flows, taken from
the operation of Internet-2 during a seven-day period. Temporal periodicities
and correlations across flows are apparent. As expected, in this case $\bZ$ can
be well approximated by a low-rank matrix, since its normalized singular values
decay rapidly to zero.}
\label{fig:fig_flows}
\end{figure}

It is not uncommon for some of the OD flow rates to experience unexpected abrupt
changes. These
so-termed \textit{traffic volume anomalies} are typically due to
(unintentional) network equipment misconfiguration or outright failure,
unforeseen behaviors following
routing policy modifications, or,
cyber attacks (e.g., DoS attacks) which aim at compromising the
services offered by the network~\cite{zggr05,lakhina,mrspmag13}.
Let $a_{f,t}$ denote the unknown amount of anomalous traffic in flow $f$ at
time $t$, which one wishes to estimate. Explicitly accounting for the presence
of anomalous flows, the measured
traffic carried by link $l$ is then given by
%
%\begin{equation}
$y_{l,t}=\sum_{f}r_{l,f}(z_{f,t}+a_{f,t}) + \epsilon_{l,t},~t=1,...,T$,
%\label{eq:y_lt}
%\end{equation}
%
where the noise variables $\epsilon_{l,t}$ capture measurement errors
and unmodeled dynamics. Traffic volume anomalies are (unsigned) sudden
changes in the traffic of OD flows, and as such their effect can span
multiple links in the network. A key difficulty in unveiling anomalies
from link-level measurements only is that  oftentimes, clearly
discernible anomalous spikes in the flow traffic can be masked through
``destructive interference'' of the superimposed OD flows~\cite{lakhina}.
An additional challenge stems from missing link-level measurements $y_{l,t}$,
an unavoidable operational reality
affecting most traffic engineering tasks that rely on (indirect) measurement
of traffic matrices~\cite{Roughan}.
To model missing link measurements, collect the tuples $(l,t)$ associated with
the
available observations $y_{l,t}$ in the set $\Omega \subseteq [1,2,...,L] \times
[1,2,...,T]$. Introducing the matrices $\bY:=[y_{l,t}],
\bE:=[\epsilon_{l,t}]\in\mathbb{R}^{L \times T}$, and $\bA:=[a_{f,t}]
\in\mathbb{R}^{F \times T}$, the (possibly incomplete) set of link-traffic
measurements
can be expressed in compact matrix form as
\begin{equation}
\cP_{\Omega}(\bY)=\cP_{\Omega}(\bX + \bR\bA + \bE)\label{eq:Y}
\end{equation}
where the sampling operator $\cP_{\Omega}(.)$ sets
the entries of its matrix argument not in $\Omega$ to zero, and keeps the rest
unchanged. Since the objective here is not
to estimate the OD flow traffic matrix $\bZ$, \eqref{eq:Y} is expressed in
terms of the nominal (anomaly-free) link-level traffic rates $\bX:=\bR\bZ$,
which inherits the low-rank property of $\bZ$. Anomalies in $\bA$ are expected
to occur sporadically over time, and last for a short time
relative to the (possibly long) measurement interval $[1,T]$.
In addition, only a small fraction of the flows is supposed to be anomalous at a
any given
time instant. This renders the anomaly matrix $\bA$
\textit{sparse} across rows (flows) and columns (time).

Recently, a natural estimator leveraging the low rank property of $\bX$ and
the sparsity of $\bA$ was put forth in~\cite{mmg13tsp}, which
can
be found at the crossroads of compressive sampling~\cite{Do06CS} and timely
low-rank plus
sparse matrix decompositions~\cite{CLMW09,CSPW11}. The idea is to fit
the incomplete data $\cP_{\Omega}(\bY)$ to the model $\bX + \bR \bA$ [cf.
\eqref{eq:Y}] in
the LS error sense, as well as minimize the
rank of $\bX$, and the number of nonzero entries of $\bA$ measured by its
$\ell_0$-(pseudo)
norm. Unfortunately, albeit natural both rank and $\ell_0$-norm criteria
are in general NP-hard to optimize.
Typically, the nuclear norm $\|\bX\|_*:=\sum_{k}\sigma_k(\bX)$ ($\sigma_k(\bX)$
denotes the $k$-th singular value of $\bX$) and the $\ell_1$-norm $\|\bA\|_1$
are adopted as surrogates~\cite{F02,CT05}, since they are the closest
\textit{convex}
approximants to $\textrm{rank}(\bX)$ and $\|\bA\|_0$, respectively.
Accordingly, one solves
\begin{equation}
\min_{\{\bX,\bA\}} \|\cP_{\Omega}(\bY - \bX -
\bR\bA)\|_{F}^{2} +\lambda_{\ast}\|\bX\|_{*} + \lambda_1\|\bA\|_1\label{eq:p1}
\end{equation}
where $\lambda_*,\lambda_1\geq 0$ are rank- and sparsity-controlling parameters.
%When an estimate $\hat{\sigma}_v^2$ of the noise variance is available,
%guidelines for selecting $\lambda_*$ and $\lambda_1$ have been proposed
%in~\cite{zlwcm10}.
While a non-smooth optimization problem, \eqref{eq:p1} is
appealing because it is convex.
An efficient accelerated proximal gradient algorithm with
quantifiable iteration complexity was developed to unveil network
anomalies~\cite{tit_recovery_2012}. Interestingly, \eqref{eq:p1}
also offers a cleansed estimate of the link-level traffic $\hat{\bX}$, that
could be subsequently utilized for network tomography tasks.
In addition, \eqref{eq:p1} \textit{jointly} exploits the spatio-temporal
correlations in link traffic as well as the sparsity of anomalies,
through an optimal single-shot estimation-detection procedure that
turns out to outperform
the algorithms  in~\cite{lakhina} and~\cite{zggr05} (the latter decouple the
estimation and
detection steps); see Fig. \ref{fig:anomography}.

\begin{figure}[ht]
\begin{minipage}[b]{0.48\linewidth}
    \centering
     \centerline{\includegraphics[width=1.05\linewidth]{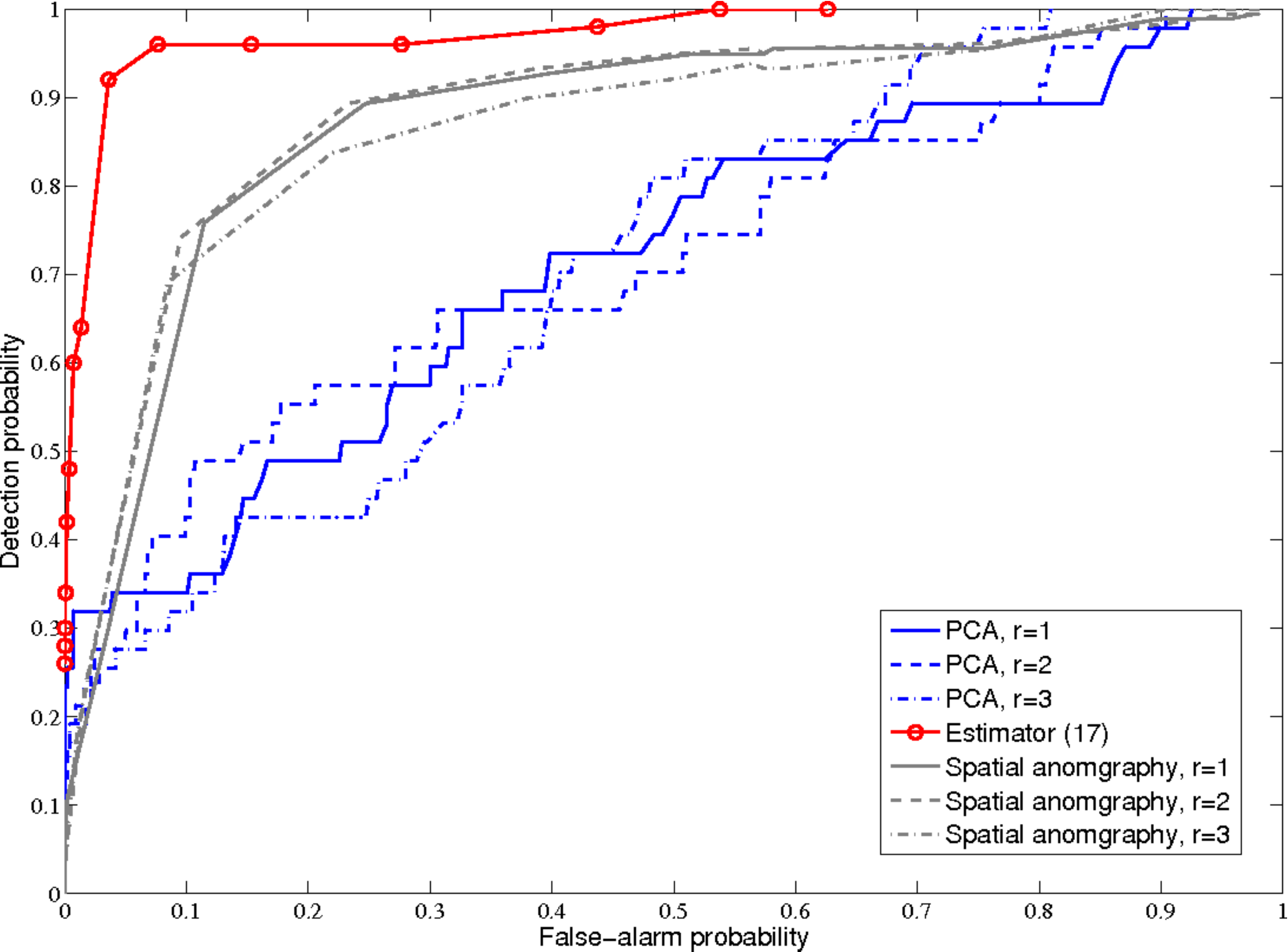}}
    \medskip
\end{minipage}
\hfill
\begin{minipage}[b]{.48\linewidth}
    \centering
     \centerline{\includegraphics[width=0.85\linewidth]{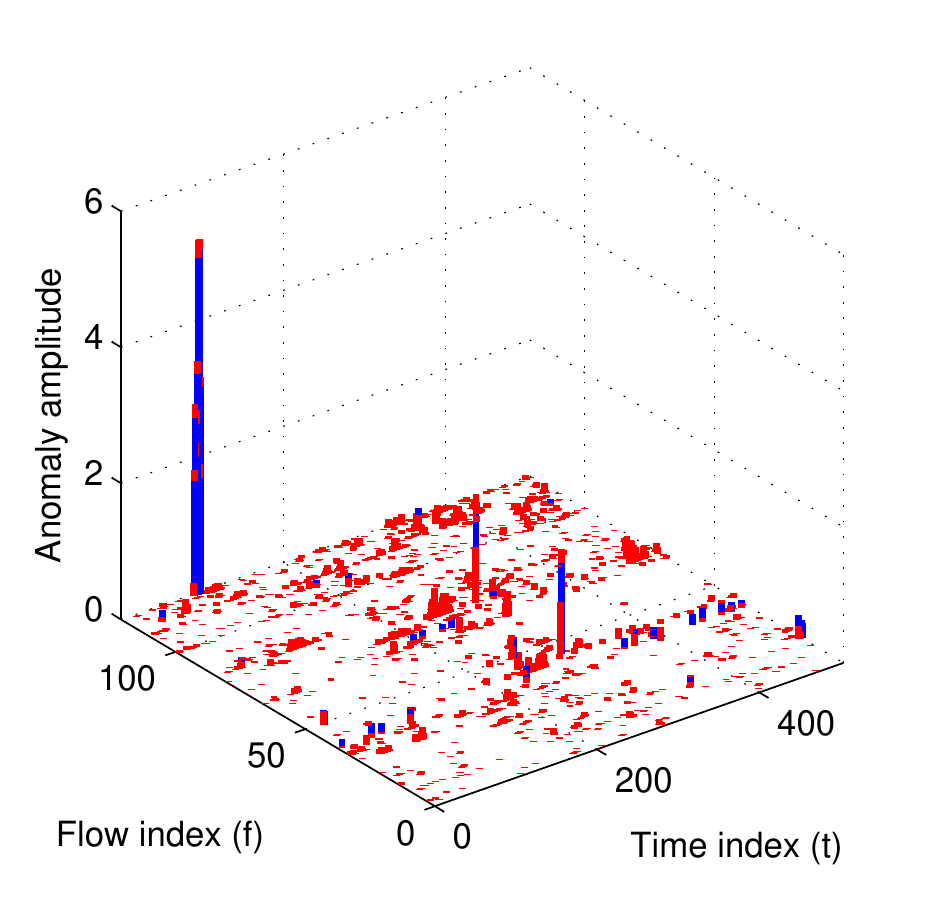}}
    \medskip
\end{minipage}

\caption{Unveiling anomalies from Internet-2 data. (Left) ROC curve comparison
    between \eqref{eq:p1} and the PCA methods in~\cite{lakhina,zggr05}, for
different values of
    the $\text{rank}(\bbZ)$. Leveraging sparsity and low rank jointly leads
    to improved performance. (Right) In red, the estimated anomaly map $\hat{\bA}$
    obtained via \eqref{eq:p1} superimposed to the ``true''
    anomalies shown in blue~\cite{mmg13jstsp}.}
\label{fig:anomography}
\end{figure}

\subsection{In-network Traffic Anomaly Detection}
\label{ssec:inp}

Implementing \eqref{eq:p1} presumes that network nodes continuously
communicate their link traffic measurements to a central monitoring station,
which uses their aggregation in $\cP_{\Omega}(\bY)$ to unveil anomalies.
While for the most part this is the prevailing operational paradigm adopted
in current networks, it is prudent to reflect on the limitations associated with this architecture.
For instance, fusing all this information may entail
excessive protocol overheads. Moreover, minimizing the exchanges of raw
measurements
may be desirable to reduce unavoidable communication errors that translate to
missing data.
Solving \eqref{eq:p1} centrally raises robustness concerns as well,
since the central monitoring station represents an isolated point of failure.

These reasons prompt one to develop \emph{fully-decentralized} iterative
algorithms for unveiling traffic anomalies, and thus embed network
anomaly detection functionality to the routers. As in Section \ref{sec:ADMM},
per iteration node $i$ carries out simple computational tasks locally,
relying on its own link count measurements (a submatrix $\bY_i$
within $\bY := [\bY_1^\top,\ldots,\bY_n^\top]^\top$ corresponding
to router $i$'s links).
Subsequently, local estimates are refined after exchanging messages
only with directly connected neighbors, which facilitates percolation
of local information to the whole network. The end goal is
for network nodes to consent on a global map of network anomalies $\hat{\bA}$,
and attain (or at least come close to) the estimation performance
of the centralized counterpart \eqref{eq:p1} which has all data
$\cP_{\Omega}(\bY)$
available.

Problem \eqref{eq:p1} is not amenable to distributed implementation
because of the non-separable nuclear norm present in the cost
function. If an upper bound $\textrm{rank}(\hat\bX)\leq
\rho$ is a priori available [recall $\hat\bX$ is the estimated link-level
traffic obtained via \eqref{eq:p1}], the search space of \eqref{eq:p1} is
effectively reduced, and one can factorize the decision variable as
$\bX=\bP\bQ^\top$, where $\bP$ and $\bQ$ are $L \times \rho$ and $T
\times \rho$ matrices, respectively. Again, it is possible to interpret
the columns of $\bX$ (viewed as points in $\mathbb{R}^L$) as belonging to
a low-rank nominal subspace, spanned by the columns of $\bP$. The rows of $\bQ$
are thus the projections of the columns of $\bX$ onto the traffic subspace.
Next, consider the following alternative characterization of
the nuclear norm (see e.g.~\cite{srebro_2005})
\begin{equation}\label{eq:nuc_nrom_def}
    \|\bX\|_*:=\min_{\{\bP,\bQ\}}~~~ \frac{1}{2}\left(\|\bP\|_F^2+\|\bQ\|_F^2
\right),\quad
    \text{s. to}~~~ \bX=\bP\bQ^\top
\end{equation}
where the optimization is over all possible bilinear factorizations
of $\bX$, so that the number of columns $\rho$ of $\bP$ and
$\mathbf{Q}$ is also a variable. Leveraging \eqref{eq:nuc_nrom_def},
the following reformulation of \eqref{eq:p1} provides an important first step
towards obtaining a decentralized algorithm for anomaly identification
\begin{equation}
    \min_{\{\bP,\bQ,\bA\}}\sum_{i=1}^n \left[ \|\cP_{\Omega_i}(\bY_i
     -\bP_i\bQ^\top-\bR_i\bA)\|_F^2+\frac{\lambda_{*}}{2n}\left(n\|\bP_i\|_F^2 +
\|\bQ\|_F^2 \right)
    + \frac{\lambda_1}{n}\|\bA\|_1\right]\label{eq:p2}
\end{equation}
which is non-convex due to the bilinear terms $\bP_i\bQ^\top$, and where
$\bR :=\left[\bR_1^\top,\ldots,\bR_n^\top\right]^\top$ is partitioned into local
routing tables available per router $i$. Adopting the separable Frobenius-norm
regularization in \eqref{eq:p2} comes with no loss of optimality relative to
\eqref{eq:p1},
provided $\textrm{rank}(\hat\bX)\leq\rho$. By finding the global minimum of
\eqref{eq:p2}
[which could entail considerably less variables than \eqref{eq:p1}],
one can recover the optimal solution of \eqref{eq:p1}. But since \eqref{eq:p2}
is
non-convex, it may have stationary points which need not be globally
optimum. As asserted in~\cite[Prop. 1]{mmg13tsp} however, if a
stationary point
$\{\bar{\bP},\bar{\bQ},\bar{\bA}\}$ of \eqref{eq:p2} satisfies
$\|\cP_{\Omega}(\bY - \bar{\bP}\bar{\bQ}^\top - \bar{\bA})\| < \lambda_*$, then
$\{\hat\bX:=\bar{\bP}\bar{\bQ}^\top,\hat{\bA}:=\bar{\bA}\}$ is the globally
optimal solution of
\eqref{eq:p1}.

To decompose the cost in \eqref{eq:p2}, in which summands inside the square
brackets
are coupled through the global variables $\{\bQ,\bA\}$, one can proceed as in
Section \ref{sec:ADMM} and
introduce auxiliary copies $\{\bQ_i,\bA_i\}_{i=1}^n$
representing local estimates of $\{\bQ,\bA\}$, one per node $i$. These
local copies along with \textit{consensus} constraints yield the decentralized
estimator
\begin{align}
    \min_{\{\bP_i,\bQ_i,\bA_i\}}& \sum_{i=1}^n \left[\|\cP_{\Omega_i}(\bY_i
    -\bP_i\bQ_i^\top-\bR_i\bA_i)\|_F^2+
    \frac{\lambda_{*}}{2n}\left(n\|\bP_i\|_F^2 + \|\bQ_i\|_F^2 \right) +
    \frac{\lambda_1}{n}\|\bA_i\|_1\right] \label{eq:p3}\\
    \text{s. to} &\quad \bQ_i=\bQ_j,\:\:\bA_i=\bA_j,\;\;
    i=1,\ldots,n,\;\;j\in{\cal  N}_i,\;\;i\neq j\nonumber
\end{align}
which follows the general form in \eqref{Eq:Separanle_Estimation}, and is
equivalent to \eqref{eq:p2} provided the
network topology graph is connected. Even though consensus
is a fortiori imposed within neighborhoods, it carries over to
the entire (connected) network and local estimates agree on the
global solution of \eqref{eq:p2}. Exploiting the separable structure
of \eqref{eq:p3} using the ADMM, a general framework for in-network
sparsity-regularized
rank minimization was put forth in~\cite{mmg13tsp}.  In a
nutshell,
local tasks per iteration $k=1,2,\ldots$ entail
solving small unconstrained quadratic programs to refine the normal subspace
$\bP_i[k]$,
in addition to soft-thresholding operations to update the anomaly maps
$\bA_i[k]$ per router. Routers exchange their estimates $\{\bQ_i[k],\bA_i[k]\}$ only
with directly connected neighbors per iteration. This way the communication
overhead remains affordable, regardless of the network size $n$.

When employed to solve non-convex problems such as \eqref{eq:p3}, so far ADMM
offers no convergence guarantees. However, there is ample
experimental evidence in the literature that supports empirical
convergence of ADMM, especially when the non-convex problem at hand exhibits
``favorable'' structure~\cite{Boyd_ADMM}. For instance, \eqref{eq:p3} is a
linearly constrained bi-convex problem with potentially good convergence
properties
-- extensive numerical tests in~\cite{mmg13tsp} demonstrate
that this is
indeed the case. While establishing convergence remains an open problem, one can
still
prove that upon convergence the distributed iterations attain
consensus and global optimality, thus offering the desirable centralized performance
guarantees~\cite{mmg13tsp}.

\subsection{RF Cartography Via Decentralized Sparse Linear Regression}
\label{ssec:rf}

In the domain of spectrum sensing for CR networks,
RF cartography amounts to constructing in a distributed fashion: i) global
power spectral density (PSD) maps capturing the distribution of radiated power
across space, time, and frequency; and ii) local channel gain (CG) maps offering the
propagation medium per frequency from each node to any point in space~\cite{cg_cartography}. These maps enable identification of
opportunistically available spectrum bands for re-use and handoff operation;
as well as localization, transmit-power estimation, and tracking of
primary user activities. While the focus here is on the
construction of PSD maps, the interested reader is referred
to~\cite{tut_rf_cartog}
for a tutorial treatment on CG cartography.

A cooperative approach to RF cartography was introduced in
\cite{bazerque}, that builds on a basis
expansion model of the PSD map $\Phi(\bbx,f)$ across space
$\bbx\in\mathbb{R}^2$, and frequency $f$.
Spatially-distributed CRs collect
smoothed periodogram samples of the received signal at given sampling
frequencies,
based on which the unknown expansion coefficients are determined.
Introducing
a virtual spatial grid of candidate source locations, the estimation task
can be cast as a linear LS problem with an augmented vector of unknown
parameters.
Still, the problem complexity (or effective degrees of freedom) can be
controlled by capitalizing
on two forms of sparsity: the first one introduced by the
narrow-band nature of transmit-PSDs relative to the broad swaths
of usable spectrum; and the second one emerging from sparsely located
active radios in the operational space (due to the grid artifact). Nonzero
entries in the parameter vector sought correspond to spatial location-frequency
band pairs
corresponding to active transmissions.
All in all, estimating the PSD map and locating the active transmitters as a
byproduct
boils down to a variable selection problem. This motivates
well employment of the ADMM and the least-absolute shrinkage and selection
operator (Lasso) for decentralized sparse linear
regression~\cite{mateos_dlasso,mmg13tsp}, an
estimator subsumed by \eqref{eq:p1}
when $\bX=\mathbf{0}_{L\times T}$, $T=1$, and matrix
$\bR$ has a specific structure that depends on the chosen bases and the path-loss
propagation
model.

Sparse total LS variants are also available to cope with
uncertainty in the regression matrix, arising due to inaccurate channel
estimation and
grid-mismatch effects~\cite{tut_rf_cartog}. Nonparametric spline-based PSD map
estimators~\cite{rf_splines} have been also shown effective in capturing general
propagation
characteristics including both shadowing and fading; see also Fig.
\ref{fig:psd_map} for
an actual PSD atlas spanning $14$ frequency sub-bands.

\begin{figure}[ht]
\centering
\centerline{\epsfig{figure=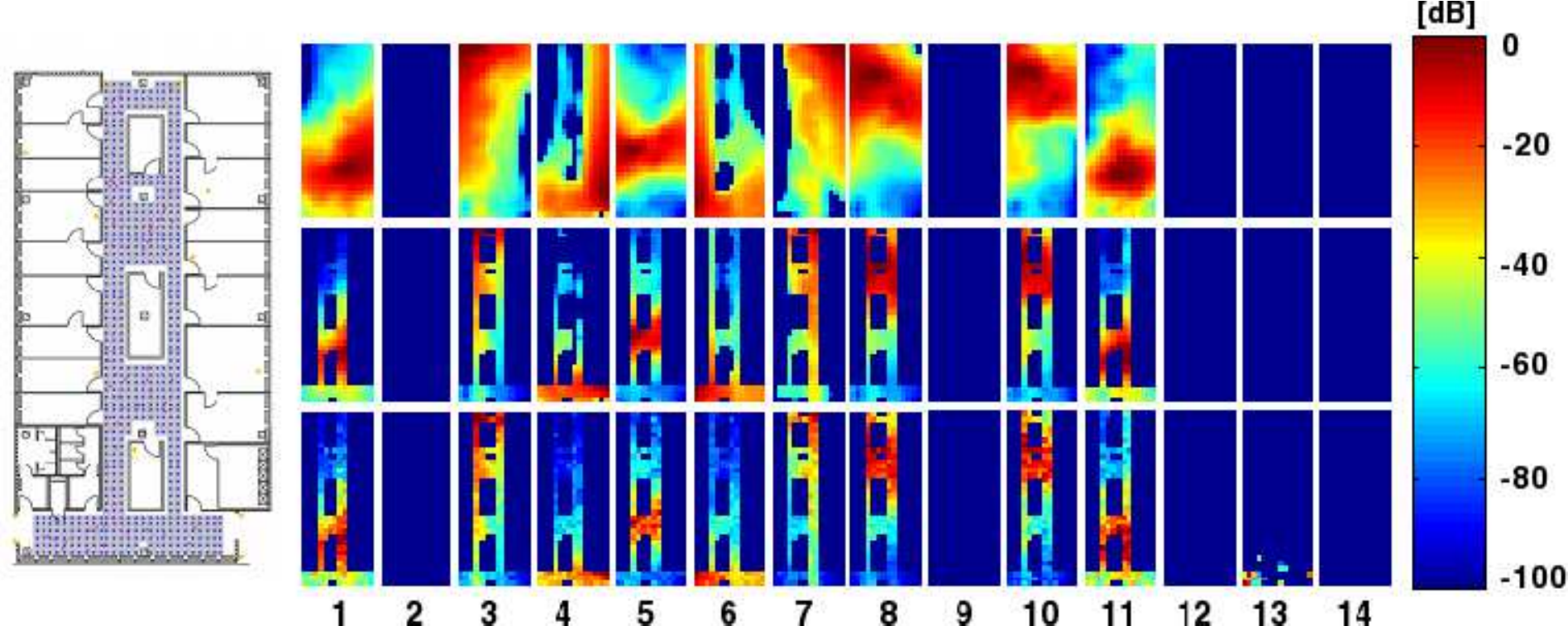,width=\linewidth}}
\caption{Spline-based RF cartography from real wireless LAN data.
    (Left) Detailed floor plan schematic including the location
    of $N=166$ sensing radios; (Right-bottom)
    original measurements spanning $14$ frequency sub-bands; (Right-center)
estimated maps over the
    surveyed area; and (Right-top) extrapolated maps. The proposed decentralized
estimator is
    capable of recovering the $9$ (out of $14$ total)
    center frequencies that are being utilized for transmission. It accurately
recovers the power
    levels in the surveyed area with a smooth extrapolation to
    zones were there are no measurements, and suggests possible
    locations for the transmitters~\cite{rf_splines}.}
\vspace{-0.8cm}  \label{fig:psd_map}
\end{figure}

\section{Convergence Analysis}
\label{sec:Convergence}

In this section we analyze the convergence and assess the rate of convergence
for the decentralized ADMM algorithm outlined in Section \ref{sec:ADMM}.
We focus on the batch learning setup, where the local cost functions are static.

\subsection{Preliminaries} \label{subsec:1.1}

\runinhead{Network model revisitied and elements of algebraic
graph theory.} Recall the network model briefly introduced  in
Section \ref{sec:Introduction}, based on a connected graph
composed of a set of $n$ nodes (agents, vertices), and a set of
$L$ edges (arcs, links).  Each edge $e= (i,j)$ represents an
ordered pair $(i,j)$ indicating that node $i$ communicates with
node $j$. Communication is assumed bidirectional so that per edge $e
= (i,j)$, the edge $e' =(j,i)$ is also present. Nodes adjacent to
$i$ are its neighbors and belong to the (neighborhood) set
$\mathcal{N}_i$. The cardinality of $|\mathcal{N}_i|$ equals the
degree $d_i$ of node $i$. Let $\bbA_s \in \mathbb{R}^{Lp \times np}$
denote the block edge source matrix, where the block
$[\bbA_s]_{e,i} = \bbI_p \in \mathbb{R}^{p\times p}$ if the edge
$e$ originates at node $i$, and is null otherwise. Likewise,
define the block edge destination matrix $\bbA_d \in
\mathbb{R}^{Lp\times np}$ where the block $[\bbA_d]_{e,j} = \bbI_p
\in \mathbb{R}^{p\times p}$ if the edge $e$ terminates at node
$j$, and is null otherwise. The so-termed extended oriented
incidence matrix can be written as $\bbE_o = \bbA_s-\bbA_d$, and
the unoriented incidence matrix as $\bbE_u = \bbA_s+\bbA_d$. The
extended oriented (signed) Laplacian is then given by $\bbL_o =
(1/2) \bbE_o^\top \bbE_o$, the unoriented (unsigned) Laplacian by
$\bbL_u = (1/2) \bbE_u^\top \bbE_u$, and the degree matrix
$\bbD=\text{diag}(d_1,\ldots,d_n)$  is
$\bbD=(1/2)(\bbL_o+\bbL_u)$. With $\Gamma_u$ denoting the
largest eigenvalue of $\bbL_u$, and $\gamma_o$ the smallest
nonzero eigenvalue of $\bbL_o$, basic results in
algebraic graph theory establish that both $\Gamma_u$ and
$\gamma_o$ are measures of network connectedness.

\runinhead{Compact learning problem representation.} With reference to the optimization
problem \eqref{Eq:Constr_Equi}, define $\bbs :=[\bbs_1^\top \ldots \bbs_n^\top]^\top \in \mathbb{R}^{np}$ concatenating all
local estimates $\bbs_i$, and $\bbz :=[\bbz_1^\top \ldots \bbz_L^\top]^\top \in \mathbb{R}^{Lp}$
concatenating all auxiliary variables $\bbz_{e}= \bbz_{i}^j$. For notational
convenience, introduce the aggregate cost
function $f: \mathbb{R}^{np} \to \mathbb{R}$ as $f(\bbs) :=
\sum_{i=1}^n f_i(\bbs_i;\bby_i)$. Using these definitions along with
the edge source and destination matrices, \eqref{Eq:Constr_Equi}
can be rewritten in compact matrix form as
\begin{equation}\label{eq:dac-mat-0}
\min_{\bbs} f(\bbs), \quad
\text{s. to }\: \bbA_s \bbs - \bbz = \mathbf{0},~\bbA_d \bbs - \bbz = \mathbf{0}. \nonumber
\end{equation}
Upon defining $\bbA := [\bbA_s^\top \bbA_d^\top]^\top \in
\mathbb{R}^{2Lp\times np}$ and  $\bbB := [-\bbI_{Lp} \:-\bbI_{Lp}]^\top$,
\eqref{eq:dac-mat-0} reduces to
\begin{equation}\label{eq:dac-mat}
    \min_{\bbs} f(\bbs), \quad
    \text{s. to } \: \bbA \bbs + \bbB\bbz = \mathbf{0}. \nonumber
\end{equation}

As in Section \ref{sec:ADMM}, consider Lagrange multipliers $\bar{\bbv}_e=\bar{\bbv}_{i}^j$ associated
with the constraints $\bbs_i = \bbs_{i}^j$, and
$\tilde{\bbv}_e=\bar{\bbv}_{i}^j$ associated with $\bbs_j = \bbs_{i}^j$. Next, define the supervectors $\bar{\bbv} :=
[\bar{\bbv}_1^\top\ldots\bar{\bbv}_L^\top]^\top \in \mathbb{R}^{Lp}$ and $\tilde{\bbv} :=
[\tilde{\bbv}_1^\top\ldots\tilde{\bbv}_L^\top]^\top \in \mathbb{R}^{Lp}$, collecting those multipliers
associated with the constraints $\bbA_s \bbs - \bbz = \mathbf{0}$ and $\bbA_d \bbs - \bbz =
\mathbf{0}$, respectively. Finally, associate multipliers $\bbv := [\bar{\bbv}^\top\:\tilde{\bbv}^\top]^\top \in
\mathbb{R}^{2Lp}$ with the constraint in \eqref{eq:dac-mat}, namely $\bbA \bbs +\bbB\bbz = \mathbf{0}$.
This way, the augmented Lagrangian function of \eqref{eq:dac-mat} is
\begin{equation} \label{eq:alf}
    L_c(\bbs,\bbz,\bbv) = f(\bbs) + \bbv^\top (\bbA\bbs+\bbB\bbz) + \frac{c}{2} \|\bbA\bbs+\bbB\bbz\|^2 \nonumber
\end{equation}
where $c > 0$ is a positive constant [cf. \eqref{Eq:Lagrangian_Function} back in Section \ref{sec:ADMM}].

\runinhead{Assumptions and scope of the convergence analysis.}
In the convergence analysis, we assume that
\eqref{Eq:Constr_Equi} has at least a pair of primal-dual
solutions. In addition, we make the following assumptions on the
local cost functions $f_i$.

% Assumption 1
\vspace{0.5em} \noindent \textbf{Assumption 1.} The local cost
functions $f_i$ are closed, proper, and convex.
%

% Assumption 2
\vspace{0.5em} \noindent \textbf{Assumption 2.} The local cost
functions $f_i$ have Lipschitz gradients, meaning there exists a
positive constant $M_{f} > 0$ such that for any node $i$ and for
any pair of points $\tilde{\bbs}_a$ and $\tilde{\bbs}_b$, it holds that
$\|\nabla f_i(\tilde{\bbs}_a) - \nabla f_i(\tilde{\bbs}_b)\| \leq
M_{f} \|\tilde{\bbs}_a - \tilde{\bbs}_b\|$.
%

% Assumption 3
\vspace{0.5em} \noindent \textbf{Assumption 3.} The local cost
functions $f_i$ are strongly convex; that is, there exists a positive
constant $m_{f} > 0$ such that for any node $i$ and for any pair
of points $\tilde{\bbs}_a$ and $\tilde{\bbs}_b$, it holds that
$(\tilde{\bbs}_a - \tilde{\bbs}_b)^\top (\nabla
f_i(\tilde{\bbs}_a) - \nabla f_i(\tilde{\bbs}_b)) \geq m_{f}
\|\tilde{\bbs}_a - \tilde{\bbs}_b\|^2$. \vspace{0.5em}

Assumption 1 implies that the aggregate function
$f(\bbs):=\sum_{i=1}^n f_i(\bbs_i;\bby_i)$ is closed, proper, and convex.
Assumption 2 ensures that the aggregate cost $f$ has Lipschitz
gradients with constant $M_f$; thus, for any pair of points $\bbs_a$ and
$\bbs_b$ it holds that
\begin{equation} \label{eq:Lipschitz_continuity}
    \|\nabla f(\bbs_a) - \nabla f(\bbs_b)\| \leq M_f \|\bbs_a - \bbs_b\|.
\end{equation}
Assumption 3 guarantees that the aggregate cost $f$ is strongly
convex with constant $m_f$; hence, for any pair of points $\bbs_a$ and $\bbs_b$
it holds that
\begin{equation} \label{eq:strong_convexity}
    \big(\bbs_a - \bbs_b\big)^\top \big( \nabla f(\bbs_a) - \nabla f(\bbs_b)\big)
    \geq m_f \|\bbs_a - \bbs_b\|^2.
\end{equation}
Observe that Assumptions 2 and 3 imply that the local cost
functions $f_i$ and the aggregate cost function $f$ are differentiable.
Assumption 1 is sufficient to prove global convergence of the
decentralized ADMM algorithm. To establish linear rate of convergence
however, one further needs Assumptions 2 and 3.

\subsection{Convergence}

In the sequel, we investigate convergence of the primal variables $\bbs(k)$ and
$\bbz(k)$ as well as the dual variable $\bbv(k)$, to their respective optimal values.
At an optimal primal solution pair $(\bbs^*,\bbz^*)$, consensus
is attained and  $\bbs^*$ is formed by $n$ stacked copies of $\tilde{\bbs}^*$,
while $\bbz^*$ also comprises $L$ stacked copies of $\tilde{\bbs}^*$, where $\tilde{\bbs}^*=\hat{\bbs}$ is an optimal solution of \eqref{Eq:Centr_Est}.
If the local cost functions are not strongly
convex, then there may exist multiple optimal primal solutions;
instead, if the local cost functions are strongly convex (i.e.,
Assumption 3 holds), the optimal primal solution is unique.

For an optimal primal solution pair $(\bbs^*,\bbz^*)$, there exist
multiple optimal Lagrange multipliers $\bbv^* := [(\bar{\bbv}^*)^\top\:
(\tilde{\bbv}^*)^\top]^\top$, where $\bar{\bbv}^* = - \tilde{\bbv}^*$~\cite{Ling2014, Shi2014}. In the following convergence analysis, we show that $\bbv(k)$
converges to one of such optimal dual solutions $\bbv^*$. In
establishing linear rate of convergence, we require that the dual
variable is initialized so that $\bbv(0)$ lies in the column
space of $\bbE_o$; and consider its convergence to a unique dual
solution $\bbv^* := [(\bar{\bbv}^*)^\top\:
(\tilde{\bbv}^*)^\top]^\top$ in which $\bar{\bbv}^*$ and
$\tilde{\bbv}^*$ also lie in the column space of $\bbE_o$. Existence and
uniqueness of such a $\bbv^*$ are also proved in~\cite{Ling2014, Shi2014}.

Throughout the analysis, define
\begin{equation}
\bbu := \left[
\begin{array}{c}
\bbs \\
\bar{\bbv}
\end{array}
\right], \quad
\bbH := \left[
\begin{array}{cc}
\frac{c}{2}\bbL_u & \mathbf{0} \\
\mathbf{0} & \frac{1}{c}\bbI_{Lp}
\end{array}
\right]. \nonumber
\end{equation}
We consider convergence of $\bbu(k)$ to its optimum
$\bbu^* :=[(\bbs^*)^\top\:(\bar{\bbv}^*)^\top]^\top$, where $(\bbs^*,\bar{\bbv}^*)$ is an optimal
primal-dual pair. The analysis is based on several contraction
inequalities, in which the distance is measured in the (pseudo)
Euclidean norm with respect to the positive semi-definite matrix
$\bbH$.

To situate the forthcoming results in context, notice that convergence of the \emph{centralized} ADMM for constrained optimization problems has been proved in e.g.,~\cite{Eckstein1992}, and its ergodic $O(1/k)$ rate of convergence is established in~\cite{He2012_SIAM, Wang2013}. For non-ergodic
convergence,~\cite{He2012} proves an $O(1/k)$ rate, and~\cite{Deng2014} improves
the rate to $o(1/k)$. Observe that in~\cite{He2012, Deng2014} the
rate refers to the speed at which the difference between two successive primal-dual iterates vanishes,
different from the speed that the primal-dual optimal iterates converge to their optima. Convergence of the decentralized ADMM is presented next in
the sense that the primal-dual iterates converge to their optima.
The analysis proceeds in four steps:

\begin{enumerate}

    \item[S1.] Show that $\|\bbu(k)-\bbu^*\|_\bbH^2$ is monotonic, namely, for
    all times $k \geq 0$ it holds that
    \begin{align}\label{eq:lemma3_contraction_main}
        \|\bbu(k+1) - \bbu^*\|_\bbH^2 \leq \|\bbu(k) - \bbu^*\|_\bbH^2 - \|\bbu(k+1) -
        \bbu(k)\|_\bbH^2.
    \end{align}

    \item[S2.] Show that $\|\bbu(k+1) - \bbu(k)\|_H^2$ is monotonically
    non-increasing, that is
    \begin{align} \label{eq:lemma4_monotonicity_main}
        \|\bbu(k+2) - \bbu(k+1)\|_\bbH^2 \leq \|\bbu(k+1) - \bbu(k)\|_\bbH^2.
    \end{align}

    \item[S3.] Derive an $O(1/k)$ rate in a non-ergodic sense based on
    \eqref{eq:lemma3_contraction_main} and
    \eqref{eq:lemma4_monotonicity_main}, i.e.,
    \begin{align}\label{eq:lemma5_sublinear_main}
        \|\bbu(k+1) - \bbu(k)\|_\bbH^2 \leq \frac{1}{k+1} \|\bbu(0) - \bbu^*\|_\bbH^2.
    \end{align}

    \item[S4.] Prove that $\bbu(k) :=[\bbs(k)^\top\:\bar{\bbv}(k)^\top]^\top$
    converges to a pair of optimal primal and dual solutions of
    (\ref{eq:dac-mat}).

\end{enumerate}

The first three steps are similar to those discussed in
\cite{He2012, Deng2014}. Proving the last step is straightforward
from the KKT conditions of (\ref{eq:dac-mat}). Under S1-S4, the main
result establishing convergence of the decentralized ADMM is as follows.

\begin{theorem} \label{theorem_1}
    If for iterations \eqref{Eq:vupdate} and \eqref{Eq:supdate}
    the initial multiplier $\bbv(0) := [\bar{\bbv}(0)^\top\:
    \tilde{\bbv}(0)^\top]^\top$ satisfies
    $\bar{\bbv}(0)=-\tilde{\bbv}(0)$, and $\bbz(0)$ is such that $\bbE_u \bbs(0) = 2\bbz(0)$, then with the ADMM penalty parameter  $c>0$ it holds under Assumption 1 that the iterates $\bbs(k)$ and $\bar{\bbv}(k)$
    converge to a pair of optimal primal and dual solutions of
    \eqref{eq:dac-mat}.
\end{theorem}

Theorem \ref{theorem_1} asserts that under proper initialization,
convergence of the decentralized ADMM only requires the local
costs $f_i$ to be closed, proper, and convex. However, it does not
specify a pair of optimal primal and dual solutions of
(\ref{eq:dac-mat}), which $(\bbs(k),\bar{\bbv}(k))$ converge to.
Indeed, $\bbs(k)$ can converge to one of the optimal primal
solutions $\bbs^*$, and $\bar{\bbv}(k)$ can converge to one of the
corresponding optimal dual solutions $\bar{\bbv}^*$. The limit
$(\bbs^*, \bar{\bbv}^*)$ is ultimately determined by the initial
$\bbs(0)$ and $\bar{\bbv}(0)$. Indeed, the conditions
in Theorem \ref{theorem_1} also guarantee ergodic and non-ergodic
$o(1/k)$ convergence rates in terms of objective error and
successive iterate differences, as proved in the recent paper
\cite{DavisYin2014}.

\subsection{Linear Rate of Convergence}

Linear rate of convergence for the \emph{centralized} ADMM is established
in~\cite{Deng2012}, and for the decentralized ADMM in~\cite{Shi2014}. Similar
to the convergence analysis of the last section, the proof includes the following steps:

\begin{enumerate}

    \item[S1'.] Show that $\|\bbu(k)-\bbu^*\|_\bbH^2$ is contractive, namely, for
    all times $k \geq 0$ it holds that
    \begin{align}\label{eq:theorem2_contraction_main}
        \|\bbu(k+1) - \bbu^*\|_\bbH^2 \leq \frac{1}{1+\delta} \|\bbu(k) - \bbu^*\|_\bbH^2
    \end{align}
    where $\delta>0$ is a constant [cf. \eqref{eqn_lemma_descent_delta}]. Note that the contraction inequality \eqref{eq:theorem2_contraction_main} implies Q-linear
    convergence of $\|\bbu(k) - \bbu^*\|_\bbH^2$.

    \item[S2'.] Show that $\|\bbs(k+1)-\bbs^*\|_\bbH^2$ is R-linearly convergent
    since it is upper-bounded by a Q-linear convergent sequence, meaning
    \begin{align}\label{eq:theorem2_bound_main}
        \|\bbs(k+1) - \bbs^*\|^2 \leq \frac{1}{m_f} \|\bbu(k) - \bbu^*\|_\bbH^2
    \end{align}
    where $m_f$ is the strong convexity constant of the aggregate cost
    function $f$.

\end{enumerate}

We now state the main result establishing linear rate of
convergence for the decentralized ADMM algorithm.

\begin{theorem} \label{theorem_2}
    If for iterations \eqref{Eq:vupdate} and \eqref{Eq:supdate} the initial multiplier $\bbv(0) := [\bar{\bbv}(0)^\top\:
    \tilde{\bbv}(0)^\top]^\top$ satisfies $\bar{\bbv}(0)=-\tilde{\bbv}(0)$; the initial auxiliary variable $\bbz(0)$ is such that $\bbE_u \bbs(0) = 2\bbz(0)$; and the initial multiplier $\bar{\bbv}(0)$ lies in the column space of $\bbE_o$,
    then with the ADMM parameter $c>0$, it holds under Assumptions 1-3, that the
    iterates $\bbs(k)$ and $\bar{\bbv}(k)$ converge R-linearly to
    $(\bbs^*,\bar{\bbv}^*)$, where $\bbs^*$ is the unique optimal primal
    solution of \eqref{eq:dac-mat}, and $\bar{\bbv}^*$ is the unique
    optimal dual solution lying in the column space of $\bbE_o$.
\end{theorem}

Theorem \ref{theorem_2} requires the local cost functions to be
closed, proper, convex, strongly convex, and have Lipschitz
gradients. In addition to the initialization dictated by Theorem
\ref{theorem_1}, Theorem \ref{theorem_2} further requires the
initial multiplier $\bar{\bbv}(0)$ to lie in the column space of
$\bbE_o$, which guarantees that $\bar{\bbv}(k)$ converges to $\bar{\bbv}^*$,
the unique optimal dual solution lying in the column space of
$\bbE_o$. The primal solution $\bbs(k)$ converges to $\bbs^*$, which is
unique since the original cost function in \eqref{Eq:Centr_Est} is
strongly convex.

Observe from the contraction inequality
(\ref{eq:theorem2_contraction_main}) that the speed of convergence
is determined by the contraction parameter $\delta$: A larger
$\delta$ means stronger contraction and hence faster convergence.
Indeed,~\cite{Shi2014} give an explicit expression of $\delta$,
that is
\begin{equation}\label{eqn_lemma_descent_delta}
    \delta = \min \left\{ \frac{(\mu-1)\gamma_o}{\mu \Gamma_u},
    \frac{2 c m_f \gamma_o}{c^2\Gamma_u \gamma_o + \mu M_f^2}
    \right\}
\end{equation}
where $m_f$ is the strong convexity constant of $f$, $M_f$ is the
Lipschitz continuity constant of $\nabla f$, $\gamma_o$ is the
smallest nonzero eigenvalue of the oriented Laplacian $\bbL_o$,
$\Gamma_u$ is the largest eigenvalue of the unoriented Laplacian
$\bbL_u$, $c$ is the ADMM penalty parameter, and $\mu>1$ is an arbitrary
constant.

As the current form of \eqref{eqn_lemma_descent_delta} does not
offer insights on how the properties of the cost functions, the
underlying network, and the ADMM parameter influence the speed of
convergence,~\cite{Ling2014,Shi2014} finds the largest value of
$\delta$ by tuning the constant $\mu$ and the ADMM
parameter $c$. Specifically,~\cite{Ling2014,Shi2014} shows that
\begin{equation*}
c = M_f \sqrt{\frac{\mu}{\Gamma_u
        \gamma_o}} \quad \text{ and }\quad
 \sqrt{\frac{1}{\mu}} =
\sqrt{\frac{1}{4}\frac{m_f^2}{M_f^2}\frac{\Gamma_L}{\gamma_L}+1}
-\frac{1}{2} \frac{m_f}{M_f} \sqrt{\frac{\Gamma_L}{\gamma_L}}
\end{equation*}
maximizes the right-hand side of \eqref{eqn_lemma_descent_delta},
so that
\begin{equation}\label{eqn_lemma_descent_comments_30}
    \delta = \frac{m_f}{M_f} \left[
    \sqrt{\frac{1}{4}\frac{m_f^2}{M_f^2}+\frac{\gamma_o}{\Gamma_u}}
    -\frac{1}{2} \frac{m_f}{M_f} \right].
\end{equation}

The best contraction parameter $\delta$ is a function of the
condition number $M_f/m_f$ of the aggregate cost function $f$, and
the condition number of the graph $\Gamma_u/\gamma_o$. Note
that we always have $\delta<1$, while small values of $\delta$
result when $M_f/m_f\gg1$ or when $\Gamma_u/\gamma_o\gg1$; that is,
when either the cost function or the graph is ill
conditioned. When the condition numbers are such that
$\Gamma_u/\gamma_o\gg M_f^2/m_f^2$, the condition number of the
graph dominates, and we obtain $\delta\approx \gamma_o/\Gamma_u$,
implying that the contraction is determined by the condition
number of the graph. When $M_f^2/m_f^2\gg \Gamma_u/\gamma_o$, the
condition number of the cost dominates and we have
$\delta \approx (m_f/M_f)\sqrt{\gamma_o/\Gamma_u}$. In the latter
case the contraction is constrained by both the condition number
of the cost function, and the condition number of the graph.

\vspace*{0.5cm}
\noindent
{\bf Acknowledgements.} The authors wish to thank the following friends,
colleagues, and co-authors who contributed to their joint publications
that the material of this chapter was extracted from: Drs. J.A. Bazerque,
A. Cano, E. Dall'Anese, S. Farahmand, N. Gatsis, P. Forero, V. Kekatos,
S.-J. Kim, M. Mardani, K. Rajawat, S. Roumeliotis, A. Ribeiro, W. Shi,
G. Wu, W. Yin, and K. Yuan. The lead author (and while with SPiNCOM all
co-authors) were supported in part from NSF grants  1202135, 1247885
1343248, 1423316, 1442686; the MURI Grant No. AFOSR FA9550-10-1-0567;
and the NIH Grant No. 1R01GM104975-01.

\end{document}